\documentclass[smallextended]{svjour3}     
\smartqed  
\usepackage{mathptmx}      
\usepackage{graphicx,epsfig}
\usepackage[all,import,rotate]{xy}
\usepackage{latexsym}
\usepackage{mathrsfs}
\usepackage{amsmath,amsfonts,amssymb,amscd,epsfig,psfrag,comment,pstricks,pst-grad,pst-plot,pst-eps,caption}

\newtheorem{Theorem}{Theorem}[section]
\newtheorem{prop}[Theorem]{Proposition}
\newtheorem{lem}[Theorem]{Lemma}

\newtheorem{cor}[Theorem]{Corollary}

\newtheorem{Specialthm}{Theorem}

\renewcommand{\theenumi}{\roman{enumi}}

\begin{document}

\title{New link invariants and Polynomials (I) \thanks{The author is supported by a grant (No. 10801021/a010402) of NSFC.}}
\subtitle{Oriented cases}


\author{Zhiqing Yang
}


\institute{Zhiqing Yang \at
              School of Mathematical Sciences, Dalian University of Technology, China \\
              Tel.: +86-411-84708351-8030\\
              \email{yangzhq@dlut.edu.cn}           
}

\date{Received: date / Accepted: date}



\maketitle

\begin{abstract}
Given any oriented link diagram, two types of new knot invariants are constructed. They satisfy some generalized skein relations. The coefficients of each invariant is from a commutative ring. Homomorphisms and representations of those rings define new link invariants. For example, the HOMFLYPT polynomial with three variables. In this sense, type one invariant is a generalization of the HOMFLYPT polynomial. Those invariants can also be modified by writhe and parameterized to get more powerful invariants. For example, the modified type one invariant distinguishes mutants, and the parameterized invariants produces information for crossing number.

\end{abstract}

\setcounter{tocdepth}{2}
\tableofcontents

\section{Introduction}\label{sec:1}

Polynomial invariants of links have a long history. In 1928, J.W. Alexander {\cite{Alexander1928}} discovered the famous Alexander polynomial. It has many connections with other topological invariants. More than 50 years later, in 1984 Vaughan Jones {\cite{Jones1987}}  discovered the Jones polynomial. Soon, the HOMFLYPT polynomial {\cite{FREYD1985}\cite{Przytycki1988}} was found. It turns out to be a generalization of both the Alexander polynomial and the Jones polynomial. There are other polynomials, for example, the Kauffman 2 variable polynomial. All those polynomials satisfy certain skein relations, which are linear equations concerning several link diagrams. A natural questions is, can they be further generalized? In this paper, we will present some new link invariants. They are natural generalizations of the HOMFLYPT polynomial, and  have 12 or 20 variables.

For simplicity, we use the following symbols to denote link diagrams.
\\
\\

\begin{figure}[h]
\begin{center}
\psset{arrowscale=2,unit=0.26}
\begin{pspicture}
(0,0)(16,0)
\psline[linewidth=1pt]{->}(0,4)(4,0)
\psline[linewidth=1pt](0,0)(1.8,1.8)
\psline[linewidth=1pt]{->}(2.2,2.2)(4,4)
\psline[linewidth=1pt]{->}(6,0)(10,4)
\psline[linewidth=1pt](6,4)(7.8,2.2)
\psline[linewidth=1pt]{->}(8.2,1.8)(10,0)
\psarc[linewidth=1pt]{<-}(14,-0.5){2}{20.0}{160.0}
\psarc[linewidth=1pt]{->}(14,4.5){2}{200.0}{340.0}
\usefont{T1}{ptm}{m}{n}
\rput(2,-1){$E_+$}
\usefont{T1}{ptm}{m}{n}
\rput(8,-1){$E_-$}
\usefont{T1}{ptm}{m}{n}
\rput(14,-1){E}
\end{pspicture}
\end{center}
\caption{The diagrams $E_+, E_-$ and $E$.}\label{f1}
\end{figure}
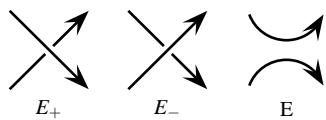

Here we use the letter $E$ because all arrows are all pointing the east direction. $+$ means positive crossing, $-$ means negative crossing. Similarly, we have the local diagrams $N_+,N_-,N,$ $W_+,W_-,W,S_+,S_-,S$. For example, the following diagrams represent $N_-, S_+$, and $S$.
\\
\\

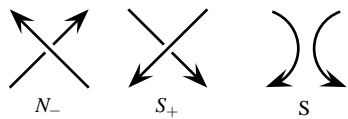
\begin{figure}[h]
\begin{center}
\psset{arrowscale=2,unit=0.26}
\begin{pspicture}
(0,0)(16,0)
\psline[linewidth=1pt]{<-}(0,4)(4,0)
\psline[linewidth=1pt](0,0)(1.8,1.8)
\psline[linewidth=1pt]{->}(2.2,2.2)(4,4)
\psline[linewidth=1pt]{<-}(6,0)(10,4)
\psline[linewidth=1pt](6,4)(7.8,2.2)
\psline[linewidth=1pt]{->}(8.2,1.8)(10,0)
\psarc[linewidth=1pt]{<-}(12.6,2){2}{-80.0}{70.0}
\psarc[linewidth=1pt]{->}(17.4,2){2}{110.0}{260.0}
\usefont{T1}{ptm}{m}{n}
\rput(2,-1){$N_-$}
\usefont{T1}{ptm}{m}{n}
\rput(8,-1){$S_+$}
\usefont{T1}{ptm}{m}{n}
\rput(14.9,-1){S}
\end{pspicture}
\end{center}
\caption{The diagrams $N_-, S_+$ and $S$.}\label{f2}
\end{figure}

Further more, we also have the followings diagrams. Here $HC$ means horizontal, and rotating clockwise. Similarly, $VT$ means vertical, and rotating anticlockwise.
\\
\\

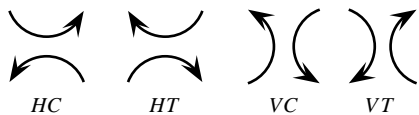
\begin{figure}[h]
\begin{center}
\psset{arrowscale=2,unit=0.26}
\begin{pspicture}
(0,0)(19,0)
\psarc[linewidth=1pt]{->}(2,-0.5){2}{20.0}{160.0}
\psarc[linewidth=1pt]{->}(2,4.5){2}{200.0}{340.0}
\psarc[linewidth=1pt]{<-}(8,-0.5){2}{20.0}{160.0}
\psarc[linewidth=1pt]{<-}(8,4.5){2}{200.0}{340.0}
\psarc[linewidth=1pt]{->}(11.5,2){2}{-70.0}{70.0}
\psarc[linewidth=1pt]{->}(16.5,2){2}{110.0}{250.0}
\psarc[linewidth=1pt]{<-}(16.6,2){2}{-70.0}{70.0}
\psarc[linewidth=1pt]{<-}(21.4,2){2}{110.0}{250.0}
\usefont{T1}{ptm}{m}{n}
\rput(2,-1){$HC$}
\usefont{T1}{ptm}{m}{n}
\rput(8,-1){$HT$}
\usefont{T1}{ptm}{m}{n}
\rput(14,-1){$VC$}
\usefont{T1}{ptm}{m}{n}
\rput(18.9,-1){$VT$}
\end{pspicture}
\end{center}
\caption{The diagrams $HC,HT,VC$ and $VT$.}\label{f3}
\end{figure}

For a local crossing $E_+$ or $E_-$ of an oriented link diagram, we propose the following new skein relations.

\noindent If the two arrows/arcs in the local diagram are from the same link component, then
$$E_{+}+bE_{-}+c_1E+c_2W+c_3HC+c_4HT+d_1VC+d_2VT=0.$$

\noindent If the two arrows/arcs are from different components, then
$$E_{+}+b'E_{-}+c_1'E+c_2'W+d'_1S+d'_2N=0.$$
We call them the {\bf type one} skein relations.

\begin{remark} For simplicity, the symbol $E_+$ ($E_-$, etc.) has many meanings in this paper. It denotes (i) the whole link diagram with the special local pattern, (ii) the local diagram contains only one crossing as in figure 1, and (iii) the value of our invariant on the diagram $E_+$. Instead of writing $$f(E_+)+bf(E_-)+c_1f(E)+c_2f(W)+c_3f(HC)+c_4f(HT)+d_1f(VC)+d_2f(VT)=0,$$ we write
   $\ \ E_++bE_-+c_1E+c_2W+c_3HC+c_4HT+d_1VC+d_2VT=0$. But sometimes, when necessary, we use $f(E_+)$ to denote the value of our invariant on the diagram $E_+$.
\end{remark}

\begin{remark}
Each diagram/term in the equations is canonically orientated as follows. For the link components/component (there are two cases) containing the arcs in the local diagram, their/its orientation is determined by the local diagrams. For all other components, the orientation is not changed. Namely, all diagrams/terms have same orientation along other components. For example, if we replace E with W, the components passing the two arrows change orientation, the other components don't. Since we distinguish the same/different component cases, there is no contradiction regarding to the orientation changes. There can be other options for orientation, please read the end of this paper.
\end{remark}

Let $A_1$ denote the commutative ring generated by $b,c_1,c_2,c_3,c_4,d_1,d_2,b',c_1',c_2',d_1',d_2'$ and $\{v_n\}_{n=1}^{\infty}$ with the following relation sets.

\noindent $R^{A_1}_1$: all generators commute.

\noindent $R^{A_1}_3$: $(1+b+d_1+d_2)v_n+(c_1+c_2+c_3+c_4)v_{n+1}=0$.

\noindent $R^{A_1}_2$:   $c_2'\overline{d}_1'=\overline{d}_1'c_2'$, $c_2'\overline{2}'=\overline{d}'_1d_2'$, $d_1'd_2'+d_2'\overline{c_1}'=d_2'd_1'+\overline{c_1}'d_2'$, $d_1'c_2'+d_2'\overline{d}_2'=\overline{1}'c_2'+d_2'c_1'$, $d_2'\overline{d}_1'=\overline{2}'c_2$, $d_2'\overline{c_2}'=\overline{c_2}'d_2'$, $c_2'\overline{c_1}'+c_1'd_2'=\overline{d}_2'd_2'+c_2'd_1'$, $c_2'\overline{d}_2'+c_1'c_2'=\overline{d}_2'c_2'+c_2'c_1', \cdots $.

The relation set $R^{A_1}_2$ is very large, and is not completely written here. The readers may find the complete description in later sections (first in the section $f_{pq}=f_{qp}$). The symbol $\overline{x}$ here denotes $b^{-1}x$, and $\overline{x}'$ denotes $b'{}^{-1}x'$. Similarly, there are commutative rings $A_2$, $A_1'$ and $A_2'$.

Here are our main theorems.

\medskip\noindent
\begin{Specialthm}[Theorem \ref{thm:main}]
{\em
For oriented link diagrams, there is a link invariant with values in $A_1$ and satisfies the following skein relations:

\noindent (1) If the two strands are from same link component, then

$E_{+}+bE_{-}+c_1E+c_2W+c_3HC+c_4HT+d_1VC+d_2VT=0.$

\noindent (2) Otherwise,
$E_{+}+b'E_{-}+c_1'E+c_2'W+d'_1S+d'_2N=0.$

The value for trivial n-component link is $v_n$.

In general, replacing $A_1$ by any homomorphic image of $A_1$, one will get a link invariant.

There is a modified invariant taking values in $A_1'$, and the value for a monotone n-component link diagram is $h(w)v_n$.}
\upshape
\end{Specialthm}
\medskip

\medskip\noindent
\begin{Specialthm}[Theorem \ref{thm:main2}]
{\em
For oriented link diagrams, there is a link invariant with values in $A_2$ and satisfies the following skein relations:

\noindent (1) If the two strands are from same link component, then

    $E_+=c_1E+c_2W+c_3HC+c_4HT+d_1VC+d_2VT$

    $E_-=\overline{c}_1E+\overline{c}_2W+\overline{c}_3HC+\overline{c}_4HT+\overline{d}_1VC+\overline{d}_2VT$

\noindent (2) Otherwise,

    $E_+=c'_1E+c'_2W+d'_1S+d'_2N$

    $E_-=\overline{c}'_1E+\overline{c}'_2W+\overline{d}'_1S+\overline{d}'_2N$

The value for a trivial n-component link diagram is $v_n$.

In general, replacing $A_2$ by any homomorphic image of $A_2$, one will get a link invariant.

There is a modified invariant taking values in $A_2'$, and the value for a trivial n-component link diagram is $v_n$.
}
\upshape
\end{Specialthm}

\begin{remark} Compare with the well-known knot polynomials, there are a few differences here. (1) The skein relation has 2 or 4 cases. (2) The coefficients now are from a commutative (or non commutative) ring, and there are some nontrivial relations among them. (3) The skein relations are more complicated. (4) The skein relation is not local here. This means for a given oriented diagram $D$, if we use the skein relation, the diagram is not only changed locally, the orientation change affects globally. To avoid contradictions, the coefficients have to satisfy certain relations. This is why we do not have a polynomial ring/invariant, but a commutative ring here.
\end{remark}

Those two types of invariants can also be modified by the writhe, like the Kauffman bracket and the Kauffman 2-variable polynomial {\cite{Kauffman1987}}. In our next paper, we shall construct similar invariants for unoriented link diagrams.

The coefficients of each invariant is from a commutative ring. Homomorphisms and representations of those rings define new link invariants. For example, if the variables in the invariants are either $0$ or invertible, one shall get knot polynomials. For example, if in the ring $A_1$ we add the following relations $c_2=c_3=c_4=d_1=d_2=c_2'=d_1'=d_2'=0$ and $b=b'$, then we get a generalized HOMFLYPT polynomial with three variables $b, c_1,c_2$. If we ask $c_1=c_1'$, then the invariant we get is equivalent to  the famous HOMFLYPT polynomial by some variable change. In this sense, it is a generalization of the HOMFLYPT polynomial. If we add other type of relations, we shall get other knot invariants. Here we list some interesting examples.

\bigskip\noindent (1) If the two strands are from the same component, we use $E_+-E_-+d_1VC+d_2VT=0$, otherwise, $E_+-E_-+d_1'S+d_2'N=0$. The relations among the coefficients are $d_1d_1'=0$, $d_1d_2'=0$, $d_2d_2'=0$, $d_1'd_2'=0$, $d_1d_1=d_2d_2=-d_1d_2$.

\noindent (2) If the two strands are from same component, we use $E_++E_-+c_1E+d_1VT=0$, otherwise, $E_++E_-+c_1'E+d_1'S=0$. The relations among the coefficients are $1d_1=1d_1'=1'd_1=d_1d_1=d_1d_1'=d_1'd_1'=0$, $(2+d_1)v_n+c_1v_{n+1}=0$.

\bigskip Those two new invariants are interesting because usually a knot polynomial does not have the $VC$ or $VT$ terms, and the second one looks similar to an ``oriented version of 2-variable Kauffman polynomial".

In fact, type one and two invariants produce many new knot polynomials, and it is a little hard to list all those polynomials.

\bigskip There are some applications of the new invariants.  For example, the modified type one invariant can distinguish mutants, which makes them more interesting. In the end, we construct new invariants with infinitely many variables (with different parameters), they are closely related to crossing number and other link invariants.

Our work was motivated by Jozef H. Przytycki and Pawel Traczyk's paper {\cite{Przytycki1988}}, and V. O. Manturov's proofs in his book {\cite{Manturov2004}}. Our construction and proof is a modification and improvement of their work.

\section{The type one invariant}\label{sec:2}

As mentioned before, we propose the following new skein relations.

\noindent If the two arrows/arcs are from the same link component, then
$$E_{+}+bE_{-}+c_1E+c_2W+c_3HC+c_4HT+d_1VC+d_2VT=0.$$

\noindent If the two arrows/arcs are from different components, then
$$E_{+}+b'E_{-}+c_1'E+c_2'W+d'_1S+d'_2N=0.$$
We call them the {\bf type one} skein relations.

\begin{remark}
If we add another variable $a$ here and use the following skein relation $$aE_{+}+bE_{-}+c_1E+c_2W+c_3HC+c_4HT+d_1VC+d_2VT=0,$$ the invariant is a little different, and possibly becomes slightly stronger. However, if we ask $a$ to have an inverse (we need only the right inverse), then we can delete this variable by changing variables, and the construction and discussions later will be much easier. Hence we do not use the variable $a$ here in this paper.
\end{remark}
\begin{remark}
There is a stronger invariant use the following skein relation:
 $$aE_{+}A+bE_{-}B+c_1EC_1+c_2WC_2+c_3HC C_3+c_4HTC_4+d_1VCD_1+d_2VTD_2=0$$
 Here $A,B,C_1,\cdots $ are new variables. However, the formulation will be much complicated. If one understand our paper, such a new construction can be similarly produced. So we do not discuss it.
\end{remark}

\begin{remark}
There is no $S$ or $N$ terms in the first equation, because if the two strands are from same component, this orientation assignment will cause contradiction in orientation. There are only four common terms between the two equations : $E_+, E_-, E$ and $W$. However, only the first three $E_+, E_-, E$ do not change the orientation of other crossings in the diagram. If one add other terms, one has to distinguish the same/different strand cases to avoid contradiction in orientation assignment.
\end{remark}

If one wants to calculate the invariant of a diagram $D$, he can start at any crossing point $p$. First, he shall determine which skein equation to use, he checks whether the two arrows/arcs of the crossing are from the same link component or not. Then, he rotates the diagram such that the crossing is either $E_+$ or $E_-$. Now he can smooth the crossing in different ways to fit in the skein equation. For example, if $p$ is a negative crossing point, and the two arcs are from the same link component, then he get: $E_{-}=-b^{-1}\{E_{+}+c_1E+c_2W+c_3HC+c_4HT+d_1VC+d_2VT\}$. Here we ask $b$ and $b'$ to have inverses. Hence if we have defined the value for $E_{+},E,\cdots $, we get the value for $E_-$. This is similar to the usually calculation of Jones polynomial by using skein relations. This also motivates us to define the invariant inductively. Such a procedure that reduces the calculation to other terms in the equation will be referred to as {\bf resolving} at $p$. We call $-b^{-1}\{E_{+}+c_1E+c_2W+c_3HC+c_4HT+d_1VC+d_2VT\}$ a {\bf linear sum}. We denote it by $f_{p}(D)$.

\subsection{$f_{pq}=f_{qp}$}\label{sec:2.1}

Given a link diagram $D$ with crossings $p_1, \cdots p_n$. Pick two crossings, say $p,q$. We can use the skein relation to resolve the diagram at a crossing $p$. The output is a linear combination of many terms. Forgetting the coefficients, each term corresponds to a link diagram $D_j$. The diagrams correspond to different ways of smoothing $p$. We denote the above by $f_p(D)=\sum \alpha_i f(D_j)$.  Each $D_j$ also has a crossing point corresponding to the crossing $q$. For each such diagram $D_j$, we resolve it at the point $q$. We shall get a new $f_q(D_j)$, a linear combination of many terms. Add the results up, we get a linear combination of linear combinations. We denote the result by $f_{pq}(D)=\sum \alpha_i f_q(D_j)$. It is the result of completely resolving at two crossing points in the order $p$ first, then $q$. Similarly, we can get another result $f_{qp}(D)$. Now, we require that if we resolve any pair $p,q$, $f_{pq}(D)=f_{qp}(D)$.

\begin{remark}
The equation $f_{pq}(D)=f_{qp}(D)$ is super important in this paper. Once this condition is satisfied, one need just a few equations to get a link invariant. We shall discuss this condition in full detail and consider several cases.
\end{remark}

\bigskip\noindent {\bf Easy cases.} $D$ is a disjoint union of two planar link diagrams $D_1$ and $D_2$, $p\in D_1$, and $q\in D_2$.

$\\$ Remember that the skein relations is as follows.

\noindent If the two arrows are from the same component, then

$E_{+}+bE_{-}+c_1E+c_2W+c_3HC+c_4HT+d_1VC+d_2VT=0$.

\noindent If the two arrows are from different components, then

$E_{+}+b'E_{-}+c_1'E+c_2'W+d'_1S+d'_2N=0$.

Hence when we apply the formula at a crossing $p$, there are two things to check, 1. the two arcs are from same/ different component, 2. the crossing is positive or negative. We call the above information the {\bf crossing pattern} of $p$.

In this case, when resolve $p$, we get diagrams $D_1, \cdots , D_k$. In all the $D_i$'s, $q$ has the same crossing pattern. In $D$, $q$ also has the same crossing pattern.

$\\$
\noindent {\bf Example 1.} If both $p,q$ are positive crossings, but for $p$, the two arrows are from same component, for $q$, the two arrows are not from same component. When we resolve $p$, we get

$-E_{+}=bE_{-}+c_1E+c_2W+c_3HC+c_4HT+d_1VC+d_2VT.$

Since we are discussing two crossings here, we use $(E,E_+)$ to denote the crossing $p$ is $E$, the crossing $q$ is $E_+$.
For each term when we resolve at $q$, we get for example:

$-b(E_{-},E_+)=bb'(E_{-},E_{-})+bc_1'(E_{-},E)+bc_2'(E_{-},W)+bd'_1(E_{-},S)+bd'_2(E_{-},N)$.

So we have the following equations.

\noindent \resizebox{12.5cm}{!} {$(E_+,E_+)=-\{b(E_{-},E_+)+c_1(E,E_+)+c_2(W,E_+)+c_3(HC,E_+)+c_4(HT,E_+) +d_1(VC,E_+)+d_2(VT,E_+)\}$}
$-b(E_{-},E_+)=bb'(E_{-},E_{-})+bc_1'(E_{-},E)+bc_2'(E_{-},W)+bd'_1(E_{-},S)+bd'_2(E_{-},N)\\$
$-c_1(E,E_+)=c_1b'(E,E_{-})+c_1c_1'(E,E)+c_1c_2'(E,W)+c_1d'_1(E,S)+c_1d'_2(E,N)\\$
$-c_2(W,E_+)=c_2b'(W,E_{-})+c_2c_1'(W,E)+c_2c_2'(W,W)+c_2d'_1(W,S)+c_2d'_2(W,N)\\$
\noindent \resizebox{12.5cm}{!} {$-c_3(HC,E_+)=c_3b'(HC,E_{-})+c_3c_1'(HC,E)+c_3c_2'(HC,W)+c_3d'_1(HC,S)+c_3d'_2(HC,N) $}
$-c_4(HT,E_+)=c_4b'(HT,E_{-})+c_4c_1'(HT,E)+c_4c_2'(HT,W)+c_4d'_1(HT,S)+c_4d'_2(HT,N)\\$
\noindent \resizebox{12.5cm}{!} {$-d_1(VC,E_+)=d_1b'(VC,E_{-})+d_1c_1'(VC,E)+d_1c_2'(VC,W)+d_1d'_1(VC,S)+d_1d'_2(VC,N) $ }
$-d_2(VT,E_+)=d_2b'(VT,E_{-})+d_2c_1'(VT,E)+d_2c_2'(VT,W)+d_2d'_1(VT,S)+d_2d'_2(VT,N)\\$

We can build a matrix for this result. We put the crossing type of the first crossing in the first column, the crossing type of the second crossing in the first row.

\begin{table}[h]
\caption{Trivial case, resolving $p$ first.}\label{tab:1}
\begin{center}
\begin{tabular}[pos]{|l|l|l|l|l|l|l|l|l|l|l| }
\hline  1 $\backslash $ 2 & $E_-$ & $E$ & $W$ & $S$ & N     \\
\hline $E_-$ & $bb'$ & $b1'$ & $b2'$  & $bd_1'$  &  $bd_2'$    \\
\hline $E$ & $1b'$ & $11'$ & $12'$  & $1d_1'$  &  $1d_2'$    \\
\hline $W$ & $2b'$ & $21'$ & $22'$  & $2d_1'$  &  $2d_2'$   \\
\hline $HC$ & $3b'$ & $31'$ & $32'$  & $3d_1'$  &  $3d_2'$   \\
\hline $HT$ & $4b'$ & $41'$ & $42'$  & $4d_1'$  &  $4d_2'$      \\
\hline $VC$ & $d_1b'$ & $d_11'$ & $d_12'$  & $d_1d_1'$  &  $d_1d_2'$    \\
\hline $VT$ & $d_2b'$ & $d_21'$ & $d_22'$  & $d_2d_1'$  &  $d_2d_2'$    \\
\hline
\end{tabular}
\end{center}
\end{table}

\begin{remark}  1. In this form/matrix, we use $1$ to denote $c_1$, similarly, $2'$ to denote $c_2'$. Later on, when the form/matrix is too wide/big, this convention makes it easier to fit an $A4$ page.

\noindent 2. Here in the form, for example the second row corresponds to $E_-$, the third row corresponds to $E$, hence the entry $b1'$ corresponds to the coefficient of $(E_-,E)$, which is $bc_1'$.
\end{remark}

Other other hand, if we resolve at $q$ first, we shall get another matrix.

\begin{table}[h]
\caption{Trivial case, resolving $q$ first.}\label{tab:2}
\begin{center}
\begin{tabular}[pos]{|l|l|l|l|l|l|l|l|l|l|l| }
\hline  1 $\backslash $ 2 & $E_-$ & $E$ & $W$ & $S$ & N     \\
\hline $E_-$ & $b'b$ & $1'b$ & $2'b$  & $d_1'b$  &  $d_2'b$    \\
\hline $E$ & $b'1$ & $1'1$ & $2'1$  & $d_1'1$  &  $d_2'1$    \\
\hline $W$ & $b'2$ & $1'2$ & $2'2$  & $d_1'2$  &  $d_2'2$   \\
\hline $HC$ & $b'3$ & $1'3$ & $2'3$  & $d_1'3$  &  $d_2'3$   \\
\hline $HT$ & $b'4$ & $1'4$ & $2'4$  & $d_1'4$  &  $d_2'4$      \\
\hline $VC$ & $b'd_1$ & $1'd_1$ & $2'd_1$  & $d_1'd_1$  &  $d_2'd_1$    \\
\hline $VT$ & $b'd_2$ & $1'd_2$ & $2'd_2$  & $d_1'd_2$  &  $d_2'd_2$    \\
\hline
\end{tabular}
\end{center}
\end{table}

Compare the results, the easiest way to make them equal is to ask the coefficients equal each other. Therefor, we ask any element from the set $\{b,c_1,c_2,c_3,c_4,d_1,d_2\}$ commutes with any element from the set $\{b',c_1',c_2',d_1',d_2'\}$.

\bigskip Let $\overline{b}=b^{-1}$, $\overline{c}_i=b^{-1}c_i$ for $i=1,2,3,4$,  $\overline{d}_i=b^{-1}d_i$ for $i=1,2$.
$\overline{b}'={b'}^{-1}$, $\overline{c}_i'={b'}^{-1}c_i'$ for $i=1,2,3,4$,  $\overline{d}_i'={b'}^{-1}d_i'$ for $i=1,2$. For any pair of such symbols $x$ and $\overline{x}$, we call them the {\bf conjugates} of each other. This has an obvious benefit as follows. In a skein relation, for example $E_{+}+bE_{-}+c_1E+c_2W+c_3HC+c_4HT+d_1VC+d_2VT=0$, we can get $E_{+}=-\{bE_{-}+c_1E+c_2W+c_3HC+c_4HT+d_1VC+d_2VT\}$ and $E_{-}=-\{\overline{b}E_{+}+\overline{c}_1E+\overline{c}_2W+\overline{c}_3HC+\overline{c}_4HT+\overline{d}_1VC+\overline{d}_2VT\}.$ This means if we change $E_+$ to $E_-$ (or $E_-$ to $E_+$), we can simply replace each $x$ to $\overline{x}$. The symmetry between them will greatly simplify our discussion later.

When we list all the subcases, we get the conclusion that any two elements from $$\{b,c_1,c_2,c_3,c_4,d_1,d_2, b',c_1',c_2',d_1',d_2'\}\cup \{\overline{b},\overline{c}_1,\overline{c}_2,\overline{c}_3,\overline{c}_4,\overline{d}_1,\overline{d}_2, \overline{b}',\overline{c}_1',\overline{c}_2',\overline{d}_1',\overline{d}_2'\}$$ are mutually commutative.

The ring here will be called the type one ring. Denote it as $A_1$.  The above is the first set of relations it satisfies. We will denote it as $R^{A_1}_1$. The subindex 1 means the first set of relations. Later on, when we build other rings, the above notations make it easier to understand the relation between the rings.

$\\$
\noindent {\bf Convention}: For convenience, in the second matrix, we exchange the order of the elements of all the terms, for example, $cd$ is changed to $dc$. So for an entry $xy$, $x$ always comes from resolving the first crossing point, $y$ always comes from resolving the second crossing point.

\bigskip
Now we are going to discuss the nontrivial cases. For simplicity, we use $a,b$ to denote the end of the first crossing $p$, and $A,B$ to denote the end of the second crossing $q$. Note that the $b$ here is not the $b$ in the skein relation. We also use them to denote the oriented strands. For example, $aAb$ means that the three arcs $a,A,b$ are from same link component, and their order is $a \to A \to b$ along the link orientation.

{}
$\\$
$\\$

\begin{figure}[h]
\begin{center}
\psset{arrowscale=2,unit=0.26}
\begin{pspicture}
(0,0)(16,0)
\psline[linewidth=1pt]{->}(0,4)(4,0)
\psline[linewidth=1pt]{->}(0,0)(4,4)
\psline[linewidth=1pt]{->}(16,0)(20,4)
\psline[linewidth=1pt]{->}(16,4)(20,0)
\usefont{T1}{ptm}{m}{n}
\rput(2,-1){The first crossing $p$}
\rput(5,4){b}
\rput(5,0.5){a}
\usefont{T1}{ptm}{m}{n}
\rput(18,-1){The second crossing $q$}
\rput(21,4){B}
\rput(21,0.5){A}
\end{pspicture}
\end{center}
\caption{The label of two crossings.}\label{f4}
\end{figure}
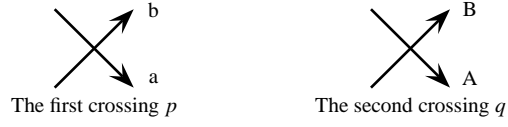

We do not distinguish positive/negative crossing type here. There is a symmetry of positive/negative crossing both in the skein relation and the diagrams. You may regard the cases we list bellow all as positive crossings. We shall tell how to deal with the other cases later.

To get all the equations $f_{pq}=f_{qp}$, we shall list all the possible cases that how the two strands of $p$ is connected to the two strands of $q$. Up to the positive/negative crossing type symmetry, there are only few interesting cases. Another observation is that we only need to discuss the cases that the four strands in the diagram above are not from 4 different components. If there are exactly two of them are from the same component, then they are not in the same local diagram. For example, if $a$ and $b$ are from same link component, then resolve $p$ would't change the crossing pattern of $q$, and vise versa. Hence the only case is (up to over/under symmetry) $aA$. Here $aA$ mean, $a,A$ are from same link component, while $b,B$ are all from different components. So there are at least three components.

If a link component contains that least two of four  ends $a,b,A,B$, we call it a {\bf complicated } component. So there is only one case when there is only one complicated component and it contains one arcs from crossings $p$ and one arcs from crossings $q$.

If there are two complicated components, then each component shall contain exactly one arc from each crossing. There is only one interesting case: case 2, $aA,bB$. The case $ab,bB$ is not interesting, it does not give any interesting equations.

If there is only one complicated component, and it contains 3 arcs from the two crossings, it may has the pattern $abA$ or $aAb$. Up to symmetry, they are the same.

If there is only one complicated component and it contains 4 arcs from the two crossings, there are only two cases up to symmetry:  $aAbB$ or $aABb$.

So, we will discuss the following five cases. 1. $aA$, 2. $aA,bB$, 3. $abA$, 4. $aAbB$, 5. $aABb$.

$\\$
\noindent {\bf Case 1, aA} If we resolve the 1st crossing point $p$ first, we shall get the followings.

$(E_+,E_+)=-\{(b'E_{-}+c_1'E+d_1'S,E_+)+(c_2'W+d_2'N,N_-)\}\\$
$-(b'E_{-},E_+)=b'\{(E_{-},b'E_{-}+c_1'E+d_1'S)+(N_+,c_2'W+d_2'N)\}\\$
$-(c'_1E,E_+)=c_1'\{(E,b'E_{-}+c_1'E+d_1'S)+(W,c_2'W+d_2'N)\}\\$
$-(d_1'S,E_+)=d_1'\{(S,b'E_{-}+c_1'E+d_1'S)+(N,c_2'W+d_2'N)\}\\$
$-(c_2'W,N_-)=c_2'\{(W,\overline{b}'N_{+}+\overline{c}_1'N+\overline{d}_2'W)+(E,\overline{c}_2'S+\overline{d}_1'E)\}\\$
$-(d_2'N,N_-)=d_2'\{(N,\overline{b}'N_{+}+\overline{c}_1'N+\overline{d}_2'W)+(S,\overline{c}_2'S+\overline{d}_1'E)\}$

\begin{table}[h]
\caption{Case aA, resolving $p$ first.}\label{tab:1a}
\begin{center}
\begin{tabular}[pos]{|l|l|l|l|l|l|l|l|l|l|l|l|}
\hline  1 $\backslash $ 2 & E & $E_-$ & $N_+$ & $N$ & S & W   \\
\hline $E$ & $2'\overline{d}_1',1'1'$ & $1'b'$ &   &  &  $2'\overline{2}',1'd_1'$ &      \\
\hline $E_-$ & $b'1'$ & $b'b'$ &   &  &  $b'd_1'$ &     \\
\hline $N_+$ &   &  &   & $b'd_2'$ &    & $b'2'$     \\
\hline $N$ &   &  & $d_2'\overline{b}'$ & $d_1'd_2',d_2'\overline{1}'$ &    & $d_1'2',d_2'\overline{d}_2'$   \\
\hline $S$ & $d_2'\overline{d}_1',d_1'1'$ & $d_1'b'$ &   &  & $d_2'\overline{2}',d_1'd_1'$ &      \\
\hline $W$ &   &   & $2'\overline{b}'$ & $2'\overline{1}',1'd_2'$ &    & $2'\overline{d}_2',1'2'$    \\
\hline
\end{tabular}
\end{center}
\end{table}

Otherwise, we shall get the followings.

$(E_+,E_+)=-\{(E_+,b'E_{-}+c_1'E+d_1'S)+(N_-,c_2'W+d_2'N)\}\\$
$-b'(E_+,E_{-})=b'\{(b'E_{-}+c_1'E+d_1'S,E_{-})+(c_2'W+d_2'N,N_+)\}\\$
$-c'_1(E_+,E)=c_1'\{(b'E_{-}+c_1'E+d_1'S,E)+(c_2'W+d_2'N,W)\}\\$
$-d_1'(E_+,S)=d_1'\{(b'E_{-}+c_1'E+d_1'S,S)+(c_2'W+d_2'N,N)\}\\$
$-c_2'(N_-,W)=c_2'\{(\overline{b}'N_{+}+\overline{c}_1'N+\overline{d}_2'W,W)+(\overline{c}_2'S+\overline{d}_1'E,E)\}\\$
$-d_2'(N_-,N)=d_2'\{(\overline{b}'N_{+}+\overline{c}_1'N+\overline{d}_2'W,N)+(\overline{c}_2'S+\overline{d}_1'E,S)\}$

\begin{table}[h]
\caption{Case aA, resolving $q$ first.}\label{tab:1b}
\begin{center}
\begin{tabular}[pos]{|l|l|l|l|l|l|l|l|l|l|l|l|}
\hline  1 $\backslash $ 2 & E & $E_-$ & $N_+$ & $N$ & S & W   \\
\hline $E$ & $\overline{d}_1'2',1'1'$ & $1'b'$ &   &  &  $\overline{d}'_1d_2',1'd_1'$ &      \\
\hline $E_-$ & $b'1'$ & $b'b'$ &   &  &  $b'd_1'$ &     \\
\hline $N_+$ &   &  &   & $\overline{b}'d_2'$ &    & $\overline{b}'2'$     \\
\hline $N$ &   &  & $d_2'b'$ & $d_2'd_1',\overline{1}'d_2'$ &    & $\overline{1}'2',d_2'1'$   \\
\hline $S$ & $\overline{2}'2,d_1'1'$ & $d_1'b'$ &   &  & $\overline{2}'d_2',d_1'd_1'$ &      \\
\hline $W$ &   &   & $2'b'$ & $\overline{d}_2'd_2',2'd_1'$ &    & $\overline{d}_2'2',2'1'$    \\
\hline
\end{tabular}
\end{center}
\end{table}

\begin{remark} Remember that in the second matrix we write the products in a new form. For example,  $(N,N_+)$ has coefficient $b'd_2'$, but we write  $d_2'b'$ in the matrix. We exchange the order of every product in this matrix so that the first symbol, for example the $d_2'$ here, is always from resolving the first crossing point $p$.
\end{remark}
\noindent The relations here are:   $c_2'\overline{d}_1'=\overline{d}_1'c_2'$, $c_2'\overline{c_2}'=\overline{d}'_1d_2'$, $b'd_2'=\overline{b}'d_2'$, $b'c_2'=\overline{b}'c_2'$, $d_2'\overline{b}'=d_2'b'$, $d_1'd_2'+d_2'\overline{c_1}'=d_2'd_1'+\overline{c_1}'d_2'$, $d_1'c_2'+d_2'\overline{d}_2'=\overline{c_1}'c_2'+d_2'c_1'$, $d_2'\overline{d}_1'=\overline{c_2}'c_2$, $d_2'\overline{c_2}'=\overline{c_2}'d_2'$, $c_2'\overline{b}'=c_2'b'$, $c_2'\overline{c_1}'+c_1'd_2'=\overline{d}_2'd_2'+c_2'd_1'$, $c_2'\overline{d}_2'+c_1'c_2'=\overline{d}_2'c_2'+c_2'c_1'$.

$\\$
\noindent {\bf Case 2, aA, bB} Resolving $p$ first, we shall get the following equations.

$(E_+,E_+)=-\{(b'E_{-}+c_1'E,E_+)+(c_2'W,W_+)+(d_2'N,N_-)+(d_1'S,S_-)\}\\$
$-(b'E_{-},E_+)=b'\{(E_{-},b'E_{-}+c_1'E)+(W_-,c_2'W)+(N_+,d_2'N)+(S_+,d_1'S)\}\\$
\noindent \resizebox{12.5cm}{!} {$-(c'_1E,E_+)=c_1'\{(E,bE_{-}+c_1E)+(W,c_2W)+(HT,c_3HC)+(HC,c_4HT)+(HC,d_1VC)+(HT,d_2VT)\}\\$}
\noindent \resizebox{12.5cm}{!} {$-(c_2'W,W_+)=c_2'\{(W,bW_{-}+c_1W)+(E,c_2E)+(HT,c_3HC)+(HC,c_4HT)+(HC,d_1VC)+(HT,d_2VT)\}\\$}
\noindent \resizebox{12.5cm}{!} {$-(d_2'N,N_-)=d_2'\{(N,\overline{b}N_{+}+\overline{c}_1N)+(S,\overline{c}_2S)+(VT,\overline{c}_3VC)+(VC,\overline{c}_4VT)+(VC,\overline{d}_1HC)+(VT,\overline{d}_2HT)\}\\$}
\noindent \resizebox{12.5cm}{!} {$-(d_1'S,S_-)=d_1'\{(S,\overline{b}S_{+}+\overline{c}_1S)+(N,\overline{c}_2N)+(VT,\overline{c}_3VC)+(VC,\overline{c}_4VT)+(VC,\overline{d}_1HC)+(VT,\overline{d}_2HT)\}$}

\begin{table}[h]
\caption{Case aA, bB, resolving $p$ first.}\label{tab:2a}
\begin{center}
\resizebox{12.5cm}{!}{\begin{tabular}[pos]{|l|l|l|l|l|l|l|l|l|l|l|l|l|l|}
\hline  1 $\backslash $ 2 & E & $E_-$ & W & $W_-$ & $N$ &$N_+$ & S & $S_+$ & HC & HT & VC & VT   \\
\hline $E$ & $2'2,1'1$ & $1'b$ &   &  &  &   & &&&&&  \\
\hline $E_-$ & $b'1'$ & $b'b'$ &   &  &  &   & &&&&&   \\
\hline $W$ &   &   & $2'1,1'2$ & $2'b$ &  &   & &&&&&   \\
\hline $W_-$ &   &   & $b'2'$ &  &  &   & &&&&&   \\
\hline $N$ &   &   &  &  & $d_1'\overline{2},d_2'\overline{1}$ & $d_2'\overline{b}$ & &&&&&   \\
\hline $N_+$ &   &   &  &  & $b'd_2'$ & & &&&&&   \\
\hline $S$ &   &   &  &  &   & & $d_1'\overline{1},d_2'\overline{2}$ & $d_1'\overline{b}$ &&&&   \\
\hline $S_+$ &   &   &  &  &   & & $b'd_1'$ &&&&&   \\
\hline $HC$ &   &   &  &  &   & &   & & & $1'4,2'4$ & $1'd_1,2'd_1$ &   \\
\hline $HT$ &   &   &  &  &   & &   & & $2'3,1'3$ & & & $2'd_2,1'd_2$ \\
\hline $VC$ &   &   &  &  &   & &   & & $d_2'\overline{d}_1,d_1'\overline{d}_1$ & & & $d_2'\overline{4},d_1'\overline{4}$  \\
\hline $VT$ &   &   &  &  &   & &   & & & $d_2'\overline{d}_2,d_1'\overline{d}_2$ & $d_1'\overline{3},d_2'\overline{3}$ &   \\
\hline
\end{tabular}}
\end{center}
\end{table}

Resolving $q$ first, we shall get the following equations.

$(E_+,E_+)=-\{(E_+,b'E_{-}+c_1'E)+(W_+,c_2'W)+(N_-,d_2'N)+(S_-,d_1'S)\}\\$
$-b'(E_+,E_{-})=b'\{(b'E_{-}+c_1'E,E_{-})+(c_2'W,W_-)+(d_2'N,N_+)+(d_1'S,S)\}\\$
\noindent \resizebox{12.5cm}{!} {$-c'_1(E_+,E)=c_1'\{(bE_{-}+c_1E,E)+(c_2W,W)+(c_3HC,HT)+(c_4HT,HC)+(d_1VC,HC)+(d_2VT,HT)\}\\$}
\noindent \resizebox{12.5cm}{!} {$-c_2'(W_+,W)=c_2'\{(bW_{-}+c_1W,W)+(c_2E,E)+(c_3HC,HT)+(c_4HT,HC)+(d_1VC,HC)+(d_2VT,HT)\}\\$}
\noindent \resizebox{12.5cm}{!} {$-d_2'(N_-,N)=d_2'\{(\overline{b}N_{+}+\overline{c}_1N,N)+(\overline{c}_2S,S)+(\overline{c}_3VC,VT)+(\overline{c}_4VT,VC)+(\overline{d}_1HC,VC)+(\overline{d}_2HT,VT)\}\\$}
\noindent \resizebox{12.5cm}{!} {$-d_1'(S_-,S)=d_1'\{(\overline{b}S_{+}+\overline{c}_1S,S)+(\overline{c}_2N,N)+(\overline{c}_3VC,VT)+(\overline{c}_4VT,VC)+(\overline{d}_1HC,VC)+(\overline{d}_2HT,VT)\}$}
\begin{table}[h]
\caption{Case aA, bB, resolving $q$ first.}\label{tab:2b}
\begin{center}
\resizebox{12.5cm}{!}{\begin{tabular}[pos]{|l|l|l|l|l|l|l|l|l|l|l|l|l|l|}
\hline  1 $\backslash $ 2 & E & $E_-$ & W & $W_-$ & $N$ &$N_+$ & S & $S_+$ & HC & HT & VC & VT   \\
\hline $E$ & $22',11'$ & $1'b'$ &   &  &  &   & &&&&&  \\
\hline $E_-$ & $b1'$ & $b'b'$ &   &  &  &   & &&&&&   \\
\hline $W$ &   &   & $12',21'$ & $2'b'$ &  &   & &&&&&   \\
\hline $W_-$ &   &   & $b2'$ &  &  &   & &&&&&   \\
\hline $N$ &   &   &  &  & $\overline{2}d_1',\overline{1}d_2'$ & $d_2'b'$ & &&&&&   \\
\hline $N_+$ &   &   &  &  & $\overline{b}d_2'$ & & &&&&&   \\
\hline $S$ &   &   &  &  &   & & $\overline{1}d_1',\overline{2}d_2'$ & $d_1'b'$ &&&&   \\
\hline $S_+$ &   &   &  &  &   & & $\overline{b}d_1'$ &&&&&   \\
\hline $HC$ &   &   &  &  &   & &   & & & $31',32'$ & $\overline{d}_1'd_2',\overline{d}_1'd_1'$ &   \\
\hline $HT$ &   &   &  &  &   & &   & & $41',42'$ & & & $\overline{d}_2'd_2',\overline{d}_2'd_1'$ \\
\hline $VC$ &   &   &  &  &   & &   & & $d_11',d_12'$ & & & $\overline{3}d_1',\overline{3}d_2'$  \\
\hline $VT$ &   &   &  &  &   & &   & & & $d_21',d_22'$ & $\overline{4}d_1',\overline{4}d_2'$ &   \\
\hline
\end{tabular}}
\end{center}
\end{table}

\noindent The relations here are:  $b'c_1'=bc_1'$, $c_2'c_1+c_1'c_2=c_1c_2'+c_2c_1'$, $b'c_2'=bc_2'$, $d_1'\overline{c_2}+d_2'\overline{c_1}=\overline{c_2}d_1'+\overline{c_1}d_2'$, $d_2'\overline{b}=d_2'b'$, $b'd_2'=\overline{b}d_2'$, $d_1'\overline{c_1}+d_2'\overline{c_2}=\overline{c_1}d_1'+\overline{c_2}d_2'$, $d_1'\overline{b}=d_1'b'$, $b'd_1'=\overline{b}d_1'$, $c_1'c_4+c_2'c_4=c_3c_1'+c_3c_2'$,$c_1'd_1+c_2'd_1=\overline{d}_1'd_2'+\overline{d}_1'd_1'$, $c_2'c_3+c_1'c_3=c_4c_1'+c_4c_2'$, $c_2'd_2+c_1'd_2=\overline{d}_2'd_2'+\overline{d}_2'd_1'$, $d_2'\overline{d}_1+d_1'\overline{d}_1=d_1c_1'+d_1c_2'$, $d_2'\overline{c_4}+d_1'\overline{c_4}=\overline{c_3}d_1'+\overline{c_3}d_2'$, $d_2'\overline{d}_2+d_1'\overline{d}_2=d_2c_1'+d_2c_2'$, $d_1'\overline{c_3}+d_2'\overline{c_3}=\overline{c_4}d_1'+\overline{c_4}d_2'$.

$\\$
\noindent {\bf Case 3, abA} Resolving $p$ first, we shall get the following equations.

$(E_+,E_+)=-\{(bE_{-}+c_1E+c_3HC+d_2VT,E_+)+(c_2W+c_4HT+d_1VC,N_-)\}\\$
$-(bE_{-},E_+)=b\{(E_{-},b'E_{-}+c_1'E+d_1'S)+(W_{-},c_2'W+d_2'N)\}\\$
$-(c_1E,E_+)=c_1\{(E,b'E_{-}+c_1'E+d_1'S)+(HT,c_2'W+d_2'N)\}\\$
$-(c_3HC,E_+)=c_3\{(HC,b'E_{-}+c_1'E+d_1'S)+(W,c_2'W+d_2'N)\}\\$
$-(d_2VT,E_+)=d_2\{(VT,b'E_{-}+c_1'E+d_1'S)+(VC,c_2'W+d_2'N)\}\\$
$-(c_2W,N_-)=c_2\{(W,\overline{b}'N_{+}+\overline{c}_1'N+\overline{d}_2'W)+(HC,\overline{c}_2'S+\overline{d}_1'E)\}\\$
$-(c_4HT,N_-)=c_4\{(HT,\overline{b}'N_{+}+\overline{c}_1'N+\overline{d}_2'W)+(E,\overline{c}_2'S+\overline{d}_1'E)\}\\$
$-(d_1VC,N_-)=d_1\{(VC,\overline{b}'N_{+}+\overline{c}_1'N+\overline{d}_2'W)+(VT,\overline{c}_2'S+\overline{d}_1'E)\}\\$

\begin{table}[h]
\caption{Case abA, resolving $p$ first.}\label{tab:3a}
\begin{center}
\begin{tabular}[pos]{|l|l|l|l|l|l|l|l|l|l|l|l|}
\hline  1 $\backslash $ 2 & E & $E_-$ & $W$ & $N$ & $N_+$ & S  \\
\hline $E$ & $4\overline{d}_1',11'$ & $1b'$ &   &  &  &  $4\overline{2}',1d_1'$ \\
\hline $E_-$ & $b1'$ & $bb'$ &   &  &   &  $bd_1'$ \\
\hline $W$ &   &   & $2\overline{d}_2',32'$ & $2\overline{1}',3d_2'$ & $2\overline{b}'$ &     \\
\hline $W_-$ &   &   & $b2'$ & $bd_2'$ &    &    \\
\hline $HC$ & $2\overline{d}_1',31'$ & $3b'$ &   &   &    & $2\overline{2}',3d_1'$     \\
\hline $HT$ &   &  & $4\overline{d}_2',12'$ & $4\overline{1}',1d_2'$ & $4\overline{b}'$ &    \\
\hline $VC$ &   &   & $d_1\overline{d}_2',d_22'$ & $d_1\overline{1}',d_2d_2'$ & $d_1\overline{b}'$ &      \\
\hline $VT$ & $d_1\overline{d}_1',d_21'$ & $d_2b'$ &   &  &  & $d_1\overline{2}',d_2d_1'$   \\
\hline
\end{tabular}
\end{center}
\end{table}

Resolving $q$ first, we shall get the following equations.

$(E_+,E_+)=-\{(E_+,b'E_{-}+c_1'E+d_1'S)+(W_+,c_2'W+d_2'N)\}\\$
$-(E_+,b'E_{-})=b'\{(bE_{-}+c_1E+c_3HC+d_2VT,E_{-})+(c_2W+c_4HT+d_1VC,N_{+})\}\\$
$-(E_+,c_1'E)=c_1'\{(bE_{-}+c_1E+c_3HC+d_2VT,E)+(c_2W+c_4HT+d_1VC,W)\}\\$
$-(E_+,d_1'S)=d_1'\{(bE_{-}+c_1E+c_3HC+d_2VT,S)+(c_2W+c_4HT+d_1VC,N)\}\\$
$-(W_+,c_2'W)=c_2'\{(bW_{-}+c_1W+c_4HT+d_1VC,W)+(c_2E+c_3HC+d_2VT,E)\}\\$
$-(W_+,d_2'N)=d_2'\{(bW_{-}+c_1W+c_4HT+d_1VC,N)+(c_2E+c_3HC+d_2VT,S)\}\\$

\begin{table}[h]
\caption{Case abA, resolving $q$ first.}\label{tab:3b}
\begin{center}
 \begin{tabular}[pos]{|l|l|l|l|l|l|l|l|l|l|l|l|}
\hline  1 $\backslash $ 2 & E & $E_-$ & $W$ & $N$ & $N_+$ & S  \\
\hline $E$ & $22',11'$ & $1b'$ &   &  &  &  $2d_2',1d_1'$ \\
\hline $E_-$ & $b1'$ & $bb'$ &   &  &   &  $bd_1'$ \\
\hline $W$ &   &   & $12',21'$ & $1d_2',2d_1'$ & $2b'$ &     \\
\hline $W_-$ &   &   & $b2'$ & $bd_2'$ &    &    \\
\hline $HC$ & $32',31'$ & $3b'$ &   &   &    & $3d_2',3d_1'$     \\
\hline $HT$ &   &  & $42',41'$ & $4d_2',4d_1'$ & $4b'$ &    \\
\hline $VC$ &   &   & $d_12',d_11'$ & $d_1d_2',d_1d_2'$ & $d_1b'$ &      \\
\hline $VT$ & $d_22',d_21'$ & $d_2b'$ &   &  &  & $d_2d_2',d_2d_1'$   \\
\hline
\end{tabular}
\end{center}
\end{table}

\noindent The relations here are:  $c_4\overline{d}_1'=c_2c_2'$, $c_4\overline{c_2}'=c_2d_2'$, $c_2\overline{d}_2'+c_3c_2'=c_1c_2'+c_2c_1'$, $c_2\overline{c_1}'+c_3d_2'=c_1d_2'+c_2d_1'$, $c_2\overline{b}'=c_2b'$, $c_2\overline{d}_1'=c_3c_2'$,  $c_2\overline{c_2}'=c_3d_2'$, $c_4\overline{d}_2'+c_1c_2'=c_4c_2'+c_4c_1'$, $c_4\overline{c_1}'+c_1d_2'=c_4d_2'+c_4d_1'$, $c_4\overline{b}'=c_4b'$, $d_1\overline{d}_2'+d_2c_2'=d_1c_2'+d_1c_1'$, $d_1\overline{c_1}'+d_2d_2'=d_1d_2'+d_1d_2'$, $d_1\overline{b}'=d_1b'$, $d_1\overline{d}_1'=d_2c_2'$, $d_1\overline{c_2}'=d_2d_2'$.

$\\$
\noindent {\bf Case 4, aAbB} Resolving $p$ first, we shall get the following equations.

\noindent \resizebox{12.5cm}{!} {$(E_+,E_+)=-\{(bE_{-}+c_1E,E_+)+(c_2W,W_+)+(c_4HT,S_-)+(c_3HC,N_-)+(d_2VT,N_-)+(d_1VC,S_-)\}\\$}
\noindent \resizebox{12.5cm}{!} {$-(bE_{-},E_+)=b\{(E_{-},bE_{-}+c_1E)+(W_-,c_2W)+(S_+,c_3HC)+(N_+,c_4HT)+(S_+,d_2VT)+(N_+,d_1VC)\}\\$}
$-(c_1E,E_+)=c_1\{(E,b'E_{-}+c_1'E)+(W,c_2'W)+(HC,d_2'N)+(HT,d_1'S)\}\\$
$-(c_2W,W_+)=c_2\{W,b'W_{-}+c_1'W)+(E,c_2'E)+(HC,d_2'N)+(HT,d_1'S)\}\\$
$-(c_4HT,S_-)=c_4\{(HT,\overline{b}'S_{+}+\overline{c}_1'S)+(HC,\overline{c}_2'N)+(E,\overline{d}_2'E)+(W,\overline{d}_1'W)\}\\$
$-(c_3HC,N_-)=c_3\{(HC,\overline{b}'N_{+}+\overline{c}_1'N)+(HT,\overline{c}_2'S)+(E,\overline{d}_1'E)+(W,\overline{d}_2'W)\}\\$
\noindent \resizebox{12.5cm}{!} {$-(d_2VT,N_-)=d_2\{(VT,\overline{b}N_{+}+\overline{c}_1N)+(VC,\overline{c}_2S)+(N,\overline{c}_3VC)+(S,\overline{c}_4VT)+(S,\overline{d}_1HC)+(N,\overline{d}_2HT)\}\\$}
\noindent \resizebox{12.5cm}{!} {$-(d_1VC,S_-)=d_1\{(VC,\overline{b}S_{+}+\overline{c}_1S)+(VT,\overline{c}_2N)+(N,\overline{c}_3VC)+(S,\overline{c}_4VT)+(S,\overline{d}_1HC)+(N,\overline{d}_2HT)\}\\$}

\begin{table}[h]
\caption{Case aAbB, resolving $p$ first.}\label{tab:4a}
\begin{center}
\resizebox{12.5cm}{!}{
\begin{tabular}[pos]{|l|l|l|l|l|l|l|l|l|l|l|l|l|l|l|l|}
\hline  1 $\backslash $ 2 & E & $E_-$ & W & $W_-$ & S & $S_+$ & N & $N_+$ & HC & HT & VC & VT  \\
\hline $E$ & $4\overline{d}_2',3\overline{d}_1'\atop 22',11'$ & $1b'$ &   &  &    &   &   &  & &  & & \\
\hline $E_-$ & $b1$ & $bb$ &  &   &  &   &   &   &   &   &   &   \\
\hline $W$ &   &   & $4\overline{d}_1',3\overline{d}_2'\atop 21',12'$  & $2b'$ &   &   &   &   &   &  &   &   \\
\hline $W_-$ &   &   & $b2$ &   &   &   &   &   &   &   &   &   \\
\hline $S$ &   &   &    &   &   &  &   &  & $d_2\overline{d}_1,d_1\overline{d}_1$ &   &   & $d_1\overline{4},d_2\overline{4}$ \\
\hline $S_+$ &   &  &  &   &  &   &   &   & $b3$ &   &   & $bd_2$ \\
\hline $N$ &  &  &  &  &   &   &  &   &   & $d_2\overline{d}_2,d_1\overline{d}_2$ & $d_1\overline{3},d_2\overline{3}$ &   \\
\hline $N_+$ &  &   &   &  &  &   &  &  &   & $b4$ & $bd_1$ &   \\
\hline $HC$ &   &  &   &   &  &  & $3\overline{1}',4\overline{2}'\atop 2d_2',1d_2'$ & $3\overline{b}'$ &   &   &   &   \\
\hline $HT$ &  &  &   &   & $4\overline{1}',3\overline{2}'\atop 2d_1',1d_1'$ & $4\overline{b}'$ &   &   &   &   &   &  \\
\hline $VC$ &  &   &    &   & $d_1\overline{1},d_2\overline{2}$ & $d_1\overline{b}$ &   &   &   &   &   &   \\
\hline $VT$ &  &   &    &   &   &  & $d_2\overline{1},d_1\overline{2}$ & $d_2\overline{b}$ &   &  &   &   \\
\hline
\end{tabular}}
\end{center}
\end{table}

Resolving $q$ first, we shall get the following equations.

\noindent \resizebox{12.5cm}{!} {$(E_+,E_+)=-\{(E_+,bE_{-}+c_1E)+(W_+,c_2W)+(S_-,c_3HC)+(N_-,c_4HT)+(S_-,d_2VT)+(N_-,d_1VC)\}\\$}
\noindent \resizebox{12.5cm}{!} {$-(E_+,bE_{-})=b\{(bE_{-}+c_1E,E_-)+(c_2W,W_-)+(c_4HT,S_+)+(c_3HC,N_+)+(d_2VT,N_{+})+(d_1VC,S_+)\}\\$}
$-(E_+,c_1E)=c_1\{(b'E_{-}+c_1'E,E)+(c_2'W,W)+(d_2'N,HT)+(d_1'S,HC)\}\\$
$-(W_+,c_2W)=c_2\{(b'W_{-}+c_1'W,W)+(c_2'E,E)+(d_1'N,HT)+(d_2'S,HC)\}\\$
$-(S_-,c_3HC)=c_3\{(\overline{b}'S_{+}+\overline{c}_1'S,HC)+(\overline{c}_2'N,HT)+(\overline{d}_2'E,E)+(\overline{d}_1'W,W)\}\\$
$-(N_-,c_4HT)=c_4\{(\overline{b}'N_{+}+\overline{c}_1'N,HT)+(\overline{c}_2'S,HC)+(\overline{d}_1'E,E)+(\overline{d}_2'W,W)\}\\$
\noindent \resizebox{12.5cm}{!} {$-(S_-,d_2VT)=d_2\{(\overline{b}S_{+}+\overline{c}_1S,VT)+(\overline{c}_2N,VC)+(\overline{c}_4VT,N)+(\overline{c}_3VC,S)+(\overline{d}_1HC,N)+(\overline{d}_2HT,S)\}\\$}
\noindent \resizebox{12.5cm}{!} {$-(N_-,d_1VC)=d_1\{(\overline{b}N_{+}+\overline{c}_1N,VC)+(\overline{c}_2S,VT)+(\overline{c}_4VT,N)+(\overline{c}_3VC,S)+(\overline{d}_1HC,N)+(\overline{d}_2HT,S)\}\\$}

\begin{table}[h]
\caption{Case aAbB, resolving $q$ first.}\label{tab:4b}
\begin{center}
\resizebox{12.5cm}{!}{\begin{tabular}[pos]{|l|l|l|l|l|l|l|l|l|l|l|l|l|l|l|l|}
\hline  1 $\backslash $ 2 & E & $E_-$ & W & $W_-$ & S & $S_+$ & N & $N_+$ & HC & HT & VC & VT  \\
\hline $E$ & $\overline{d}_2'3,\overline{d}_1'4 \atop 2'2,1'1$ & $1b$ &   &  &    &   &   &  & &  & & \\
\hline $E_-$ & $b'1$ & $bb$  &  &   &  &   &   &   &   &   &   &   \\
\hline $W$ &   &   & $\overline{d}_1'3,\overline{d}_2'4 \atop 2'1,1'2$  & $2b$ &   &   &  &   &   &   &  &   \\
\hline $W_-$ &  &  & $b'2$  &   &   &  &  &  &   &   &   &  \\
\hline $S$ &  &   &   &  &   &   &   &   & $\overline{1}'3,\overline{2}'4\atop d_2'2,d_1'1$ &   &   & $\overline{1}d_2,\overline{2}d_1$ \\
\hline $S_+$ &   &   &   &   &   &   &   &   & $\overline{b}'3$ &  &   & $\overline{b}d_2$ \\
\hline $N$ &  &   &   &  &  &   &   &   &   & $\overline{1}'4,\overline{2}'3\atop d_1'2,d_2'1$ & $\overline{1}d_1,\overline{2}d_2$ &   \\
\hline $N_+$ &  &   &   &  &  &   &   &   &   &  $\overline{b}'4$ & $\overline{b}d_1$ &  \\
\hline $HC$ &   &   &    &   &   &   & $\overline{d}_1d_2,\overline{d}_1d_1$ & $3b$ &   &   &   &   \\
\hline $HT$ &   &   &    &   & $\overline{d}_2d_2,\overline{d}_2d_1$ & $3b$ &   &   &   &   &   &  \\
\hline $VC$ &   &   &    &   & $\overline{3}d_1,\overline{3}d_2$ & $d_1b$ &   &   &   &   &   &  \\
\hline $VT$ &   &   &   &   &   &   & $\overline{4}d_1,\overline{4}d_2$ & $d_2b$ &   &   &   &  \\
\hline
\end{tabular}}
\end{center}
\end{table}

\noindent The relations here are: $c_4\overline{d}_2'+c_3\overline{d}_1'+ c_2c_2'+c_1c_1'=\overline{d}_2'c_3+\overline{d}_1'c_4 + c_2'c_2+c_1'c_1$, $c_1b'=c_1b$, $bc_1=b'c_1$, $c_4\overline{d}_1'+c_3\overline{d}_2'+ c_2c_1'+c_1c_2'=\overline{d}_1'c_3+\overline{d}_2'c_4 + c_2'c_1+c_1'c_2$, $c_2b'=c_2b$, $bc_2=b'c_2$, $d_2\overline{d}_1+d_1\overline{d}_1=\overline{1}'3+\overline{2}'4+ d_2'2+d_1'1$, $d_1\overline{4}+d_2\overline{4}=\overline{1}d_2+\overline{2}d_1$,
$bc_3=\overline{b}'c_3$, $c_bd_2=\overline{b}d_2$, $d_2\overline{d}_2+d_1\overline{d}_2=\overline{1}'c_4+\overline{2}'c_3+ d_1'c_2+d_2'c_1$, $d_1\overline{c_3}+d_2\overline{c_3}= \overline{c_1}d_1+\overline{c_2}d_2$, $bc_4= \overline{b}'c_4$, $bd_1=\overline{b}d_1$, $c_3\overline{1}'+c_4\overline{2}'+ c_2d_2'+c_1d_2'=\overline{d}_1d_2+\overline{d}_1d_1$, $c_4\overline{b}'=c_3b$, $d_1\overline{c_1}+d_2\overline{c_2}=\overline{c_3}d_1+\overline{c_3}d_2$, $d_1\overline{b}=d_1b$, $d_2\overline{c_1}+d_1\overline{c_2}=\overline{c_4}d_1+\overline{c_4}d_2$, $d_2\overline{b}=d_2b$.

$\\$
\noindent {\bf Case 5, aABb} Resolving $p$ first, we shall get the following equations.

$(E_+,E_+)=-\{(bE_{-}+c_1E+c_4HT+d_1VC,E_+)+(c_2W+c_3HC+d_2VT,W_+)\}\\$
$-(bE_{-},E_+)=b\{(E_{-},bE_{-}+c_1E+c_3HC+d_2VT)+(W_{-},c_2W+c_4HT+d_1VC)\}\\$
$-(c_1E,E_+)=c_1\{(E,bE_{-}+c_1E+c_3HC+d_2VT)+(HC,c_2W+c_4HT+d_1VC)\}\\$
$-(c_4HT,E_+)=c_4\{(HT,bE_{-}+c_1E+c_3HC+d_2VT)+(W,c_2W+c_4HT+d_1VC)\}\\$
$-(d_1VC,E_+)=d_1\{(VC,bE_{-}+c_1E+c_3HC+d_2VT)+(VT,c_2W+c_4HT+d_1VC)\}\\$
$-(c_2W,W_+)=c_2\{(W,bW_{-}+c_1W+c_4HT+d_1VC)+(HT,c_2E+c_3HC+d_2VT)\}\\$
$-(c_3HC,W_+)=c_3\{(HC,bW_{-}+c_1W+c_4HT+d_1VC)+(E,c_2E+c_3HC+d_2VT)\}\\$
$-(d_2VT,W_+)=d_2\{(VT,bW_{-}+c_1W+c_4HT+d_1VC)+(VC,c_2E+c_3HC+d_2VT)\}$

\begin{table}[h]
\caption{Case aABb, resolving $p$ first.}\label{tab:5a}
\begin{center}
\begin{tabular}[pos]{|l|l|l|l|l|l|l|l|l|}
\hline  1 $\backslash $ 2 & $E_-$ & E & $W_-$ & W & HC & HT & VC & VT \\
\hline $E_-$ & $bb$ & $b1$ &   &  & $b3$ &  &  & $bd_2$  \\
\hline E & $1b$ & $32,11$ &   &  &  $33,13$ &&& $3d_2,1d_2$    \\
\hline $W_-$ &  &   &   & $b2$ &  & $b4$ & $bd_1$ &    \\
\hline W &   &   & $2b$ & $21,42$ &  & $24,43$ & $2d_1,4d_1$ &    \\
\hline HC &   &   & $3b$ & $31,12$  &    & $34,14$ & $3d_1,1d_1$ &     \\
\hline HT & $4b$  & $22,41$ &   &  & $23,43$ &&& $2d_2,4d_2$  \\
\hline VC & $d_1b$ & $d_11,d_22$ &   &  & $d_13,d_23$ &&& $d_1d_2,d_2d_2$  \\
\hline VT &   &   & $d_2b$ & $d_21,d_12$ &  & $d_14,d_24$ & $d_1d_1,d_2d_1$ &     \\
 \hline
\end{tabular}
\end{center}
\end{table}
Resolving $q$ first, we shall get the following equations.

$(E_+,E_+)=-\{(E_+,bE_{-}+c_1E+c_3HC+d_2VT)+(W_+,c_2W+c_4HT+d_1VC)\} \\$
$-(E_+,bE_{-})=b\{(bE_{-}+c_1E+c_4HT+d_1VC,E_{-})+(c_2W+c_3HC+d_2VT,W_{-})\}\\$
$-(E_+,c_1E)=c_1\{(bE_{-}+c_1E+c_4HT+d_1VC,E)+(c_2W+c_3HC+d_2VT,HT)\}\\$
$-(E_+,c_3HC)=c_3\{(bE_{-}+c_1E+c_4HT+d_1VC,HC)+(c_2W+c_3HC+d_2VT,W)\}\\$
$-(E_+,d_2VT)=d_2\{(bE_{-}+c_1E+c_4HT+d_1VC,VT)+(c_2W+c_3HC+d_2VT,VC)\}\\$
$-(W_+,c_2W)=c_2\{(bW_{-}+c_1W+c_3HC+d_2VT,W)+(c_2E+c_4HT+d_1VC,HC)\}\\$
$-(W_+,c_4HT)=c_4\{(bW_{-}+c_1W+c_3HC+d_2VT,HT)+(c_2E+c_4HT+d_1VC,E)\}\\$
$-(W_+,d_1VC)=d_1\{(bW_{-}+c_1W+c_3HC+d_2VT,VC)+(c_2E+c_4HT+d_1VC,VT)\}.\\$

\begin{table}[h]
\caption{Case aABb, resolving $q$ first.}\label{tab:5b}
\begin{center}
\begin{tabular}[pos]{|l|l|l|l|l|l|l|l|l|}
\hline  1 $\backslash $ 2 & $E_-$ & E & $W_-$ & W & HC & HT & VC & VT \\
\hline $E_-$ & $bb$ & $b1$ &   &  & $b3$ &  &  & $bd_2$  \\
\hline E & $1b$ & $42,11$ &   &  &  $22,13$ &&& $2d_1,1d_2$    \\
\hline $W_-$ &  &   &   & $b2$ &  & $b3$ & $bd_1$ &    \\
\hline W &   &   & $2b$ & $12,23$ &  & $13,21$ & $1d_1,2d_2$ &    \\
\hline HC &   &   & $3b$ & $32,33$  &    & $33,31$ & $3d_1,4d_2$ &     \\
\hline HT & $4b$  & $43,41$ &   &  & $42,43$ &&& $4d_2,4d_1$  \\
\hline VC & $d_1b$ & $d_13,d_11$ &   &  & $d_12,d_13$ &&& $d_1d_1,d_1d_2$  \\
\hline VT &   &   & $d_2b$ & $d_22,d_23$ &  & $d_23,d_21$ & $d_2d_1,d_2d_2$ &     \\
 \hline
\end{tabular}
\end{center}
\end{table}

\noindent The relations here are: $c_3c_2=c_4c_2$, $c_3c_3=c_2c_2$, $c_3d_2=c_2d_1$, $bc_3=bc_4$, $c_2c_1+c_4c_2=c_1c_2+c_2c_3$, $c_2c_4+c_4c_3=c_1c_3+c_2c_1$, $c_2d_1+c_4d_1=c_1d_1+c_2d_2$, $c_3c_1+c_1c_2=c_3c_2+c_3c_3$, $c_3c_4+c_1c_4=c_3c_4+c_3c_1$, $c_1d_1=c_4d_2$, $c_4c_3=c_2c_2$, $c_2c_3=c_4c_2$, $c_2d_2=c_4d_1$, $d_2c_2=d_1c_3$, $d_2c_3=d_1c_2$, $d_2d_2=d_1d_1$, $d_2c_1+d_1c_2=d_2c_2+d_2c_3$, $d_1c_4+d_2c_4=d_2c_3+d_2c_1$.

\bigskip In short, here are all the relations if the two crossings are all positive.

 \noindent {\bf Case 1:}  $c_2'\overline{d}_1'=\overline{d}_1'c_2'$, $c_2'\overline{c_2}'=\overline{d}'_1d_2'$, $b'd_2'=\overline{b}'d_2'$, $b'c_2'=\overline{b}'c_2'$, $d_2'\overline{b}'=d_2'b'$, $d_1'd_2'+d_2'\overline{c_1}'=d_2'd_1'+\overline{c_1}'d_2'$, $d_1'c_2'+d_2'\overline{d}_2'=\overline{c_1}'c_2'+d_2'c_1'$, $d_2'\overline{d}_1'=\overline{c_2}'c_2$, $d_2'\overline{c_2}'=\overline{c_2}'d_2'$, $c_2'\overline{b}'=c_2'b'$, $c_2'\overline{c_1}'+c_1'd_2'=\overline{d}_2'd_2'+c_2'd_1'$, $c_2'\overline{d}_2'+c_1'c_2'=\overline{d}_2'c_2'+c_2'c_1'$.

 \noindent {\bf Case 2:}  $b'c_1'=bc_1'$, $c_2'c_1+c_1'c_2=c_1c_2'+c_2c_1'$, $b'c_2'=bc_2'$, $d_1'\overline{c_2}+d_2'\overline{c_1}=\overline{c_2}d_1'+\overline{c_1}d_2'$, $d_2'\overline{b}=d_2'b'$, $b'd_2'=\overline{b}d_2'$, $d_1'\overline{c_1}+d_2'\overline{c_2}=\overline{c_1}d_1'+\overline{c_2}d_2'$, $d_1'\overline{b}=d_1'b'$, $b'd_1'=\overline{b}d_1'$, $c_1'c_4+c_2'c_4=c_3c_1'+c_3c_2'$,$c_1'd_1+c_2'd_1=\overline{d}_1'd_2'+\overline{d}_1'd_1'$, $c_2'c_3+c_1'c_3=c_4c_1'+c_4c_2'$, $c_2'd_2+c_1'd_2=\overline{d}_2'd_2'+\overline{d}_2'd_1'$, $d_2'\overline{d}_1+d_1'\overline{d}_1=d_1c_1'+d_1c_2'$, $d_2'\overline{c_4}+d_1'\overline{c_4}=\overline{c_3}d_1'+\overline{c_3}d_2'$, $d_2'\overline{d}_2+d_1'\overline{d}_2=d_2c_1'+d_2c_2'$, $d_1'\overline{c_3}+d_2'\overline{c_3}=\overline{c_4}d_1'+\overline{c_4}d_2'$.

 \noindent {\bf Case 3:}  $c_4\overline{d}_1'=c_2c_2'$, $c_4\overline{c_2}'=c_2d_2'$, $c_2\overline{d}_2'+c_3c_2'=c_1c_2'+c_2c_1'$, $c_2\overline{c_1}'+c_3d_2'=c_1d_2'+c_2d_1'$, $c_2\overline{b}'=c_2b'$, $c_2\overline{d}_1'=c_3c_2'$,  $c_2\overline{c_2}'=c_3d_2'$, $c_4\overline{d}_2'+c_1c_2'=c_4c_2'+c_4c_1'$, $c_4\overline{c_1}'+c_1d_2'=c_4d_2'+c_4d_1'$, $c_4\overline{b}'=c_4b'$, $d_1\overline{d}_2'+d_2c_2'=d_1c_2'+d_1c_1'$, $d_1\overline{c_1}'+d_2d_2'=d_1d_2'+d_1d_2'$, $d_1\overline{b}'=d_1b'$, $d_1\overline{d}_1'=d_2c_2'$, $d_1\overline{c_2}'=d_2d_2'$.

 \noindent {\bf Case 4:}  $c_4\overline{d}_2'+c_3\overline{d}_1'+ c_2c_2'+c_1c_1'=\overline{d}_2'c_3+\overline{d}_1'c_4 + c_2'c_2+c_1'c_1$, $c_1b'=c_1b$, $bc_1=b'c_1$, $c_4\overline{d}_1'+c_3\overline{d}_2'+ c_2c_1'+c_1c_2'=\overline{d}_1'c_3+\overline{d}_2'c_4 + c_2'c_1+c_1'c_2$, $c_2b'=c_2b$, $bc_2=b'c_2$, $d_2\overline{d}_1+d_1\overline{d}_1=\overline{1}'3+\overline{2}'4+ d_2'2+d_1'1$, $d_1\overline{4}+d_2\overline{4}=\overline{1}d_2+\overline{2}d_1$,
$bc_3=\overline{b}'c_3$, $c_bd_2=\overline{b}d_2$, $d_2\overline{d}_2+d_1\overline{d}_2=\overline{1}'c_4+\overline{2}'c_3+ d_1'c_2+d_2'c_1$, $d_1\overline{c_3}+d_2\overline{c_3}= \overline{c_1}d_1+\overline{c_2}d_2$, $bc_4= \overline{b}'c_4$, $bd_1=\overline{b}d_1$, $c_3\overline{1}'+c_4\overline{2}'+ c_2d_2'+c_1d_2'=\overline{d}_1d_2+\overline{d}_1d_1$, $c_4\overline{b}'=c_3b$, $d_1\overline{c_1}+d_2\overline{c_2}=\overline{c_3}d_1+\overline{c_3}d_2$, $d_1\overline{b}=d_1b$, $d_2\overline{c_1}+d_1\overline{c_2}=\overline{c_4}d_1+\overline{c_4}d_2$, $d_2\overline{b}=d_2b$.

 \noindent {\bf Case 5:}  $c_3c_2=c_4c_2$, $c_3c_3=c_2c_2$, $c_3d_2=c_2d_1$, $bc_3=bc_4$, $c_2c_1+c_4c_2=c_1c_2+c_2c_3$, $c_2c_4+c_4c_3=c_1c_3+c_2c_1$, $c_2d_1+c_4d_1=c_1d_1+c_2d_2$, $c_3c_1+c_1c_2=c_3c_2+c_3c_3$, $c_3c_4+c_1c_4=c_3c_4+c_3c_1$, $c_1d_1=c_4d_2$, $c_4c_3=c_2c_2$, $c_2c_3=c_4c_2$, $c_2d_2=c_4d_1$, $d_2c_2=d_1c_3$, $d_2c_3=d_1c_2$, $d_2d_2=d_1d_1$, $d_2c_1+d_1c_2=d_2c_2+d_2c_3$, $d_1c_4+d_2c_4=d_2c_3+d_2c_1$.

\begin{remark}
We list here the nontrivial relations when the two crossings are all positive. How to handle the negative crossings? Well, this is very simple. For example, when the first crossing is changed to negative, we change the corresponding coefficient $x$ to $\overline{x}$. Then in the matrices we get, we change each entry. For example, $xy$ is changed to $\overline{x}y$. If the second crossing is changed to negative, we change the second symbol. For example,  $d_2\overline{b}=d_2c_1$ is changed to $d_2b=d_2\overline{c}_1$. If both the crossings are negative, we change both the symbols. In this sense, we say the relation is closed under conjugation. This means whenever we have a relation $xy=zw$, we then always have $x\overline{y}=z\overline{w}, \overline{x}y=\overline{z}w, \overline{xy}=\overline{zw}$. The collection of all nontrivial relations above and their conjugates will be denoted by $R^{A_1}_2$.
\end{remark}

\subsection{The construction and proofs for the type one invariant}
To define the invariant on any oriented link diagram $D$, we shall first assume/add some additional data.

   (1) Suppose each link component has an {\bf orientation}. This is already given.

   (2) {\bf Order} the link components by integers: 1,2, $\cdots$, m.

   (3) On each component $k_i$, pick a {\bf base point} $p_i$.

An oriented link diagram with ordering of link components and bas points is called a {\bf marked diagram}. Now, we go through component $k_1$ from $p_1$ along its orientation. When we finish $k_1$, we shall pass to $k_2$ starting from $p_2$, $\cdots$.

\begin{definition} A crossing point is called {\bf bad} if it is first passed over, otherwise, it is called {\bf good}. A link diagram contains only good crossings is called a monotone or ascending diagram.
\end{definition}

Given a monotone diagram, each link component $k_i$ can be regarded as a map $k_i: S^1 \to R^2 \times R$, and the $S^1$ can be divided into two arcs $\alpha \cup \beta$, such that, (1) the map $\beta \to R^2 \times R \to R^2$ is an immersion, (2) different points in $\beta$ has different $R$ coordinates (the third coordinate in $R^2 \times R =R^3$), hence $\beta \to R^2 \times R \to R$ is monotonously increasing, (3) the image of $\alpha$ is vertical, i.e. its projection on $R^2$ is one single point, (4) any point in $k_i$ has smaller $R$ coordinate than the points in $k_{i+1}$. The set of maps $\{k_i\}$ is called a {\bf geometric realization} of a monotone diagram.

\begin{lemma} A monotone diagram corresponds to a trivial link.
\end{lemma}

\noindent We do not use this lemma in this paper. It will help the readers to understand why we define the value for monotone diagram to be $v_n$. The proof is easy. We leave it as an exercise.

\bigskip Now we are going to construct a new link invariant for oriented link diagrams. We shall prove it is an invariant at two levels, first, at diagram level, second, at link invariant level. In this paper, a link diagram invariant mean this invariant is well-defined for a fixed diagram, no matter which crossing point you resolve first, but it may change under Reidemeister moves. A link invariant means a link diagram invariant which is also invariant under Reidemeister moves. We shall use the following proposition. The statement (1) in this proposition and its corollary 2.15 tells us that this is an invariant at diagram level. The whole proposition tells us that this is a link invariant.

For a given marked link diagram, we can define an ordered pair $(c,d)$ of integers, called its index. Here $c$ is the crossing number of the diagram, and $d$ is the number of bad points of the diagram. $(c,d)<(c',d')$ if $c<c'$, or $c=c'$ and $d<d'$. Let $S(c,d)$ denote the set of all marked link diagrams with indices $\leq (c,d)$. Note that $S(c,0)$ contains exactly the monotone diagrams with $n$ crossing points. To prove the proposition, let's first study the skein relations. Take $f(E_{+})+bf(E_{-})+c_1f(E)+c_2f(W)+c_3f(HC)+c_4f(HT)+d_1f(VC)+d_2f(VT)=0$ for example, each term has a link diagram corresponding to it. If the diagram $E_+$ is marked, then the link diagram $E_{-}$ is canonically marked using the same orientation, order, base points as $E_+$. The other diagrams are canonically orientated only. Suppose the marked link diagram $E_{+}$ has index $(c,d)$, then $E_{-}$ has index $(c,d+1)$ or $(c,d-1)$, and all other diagrams has crossing number $c-1$. As we will show later, the invariant actually does not depend on the order and base points of the link diagram, this tells us that we can construct the invariant and prove its properties use induction on the index pair $(c,d)$. For example, if $E_{-}$ has index $(c,d+1)$, and the invariant is already defined for any diagram with index $\leq (c,d)$, then $f(E_+)$ is canonically defined, and all other diagrams in the equation is oriented  and with one fewer crossing, hence the values are also uniquely defined. Then the skein relation defines $f(E_-)$ uniquely. We shall use this as the definition of $f(E_-)$. At the same time, we say that if we resolve at this bad point of $E_-$, the skein relation is satisfied.

\begin{prop} There is an invariant defined on marked link diagrams, satisfies the following properties.

\noindent (0) The value for any marked link diagram is uniquely defined.

\noindent (1) Resolving at any bad point, the invariant satisfies the type one skein relations.

\noindent (2) It is invariant under base point changes.

\noindent (3) It is invariant under Reidemeister moves.

\noindent (4) It is invariant under changing order of components.

\end{prop}

\begin{remark}
We shall prove it inductively and use the following notations.

\noindent $D_{c,d}$ means that the value for any marked link diagram in $S(c,d)$ is uniquely defined.

\noindent $Sf_{c,d}$ means that for any given marked diagram $E_{+}$ (or $E_-$) in $S(c,d)$, if $E_{-}$ (or $E_+$) has lower index, then the skein relations is satisfied when resolving at the first bad point.

\noindent $Sb_{c,d}$ means that for any given marked diagram $E_{+}$ (or $E_-$) with $\leq c$ crossings, then the skein relations is satisfied when resolving at any bad point.

\noindent $S_c$ means that for any given marked diagram with $\leq c$ crossings, then the skein relations is satisfied when resolving at any crossing point.

\noindent $B_{c,d}$ means that the value is invariant under base point changes for any two marked link diagrams in $S(c,d)$.

\noindent $R^i_{c,d}$, $i=1,2,3$, means that for any two marked diagrams $D, D'$ in $S(c,d)$ which are connected by a Reidemeister-$i$ move, $f(D)=f(D')$. Here, two marked diagrams $D, D'$ in $S(c,d)$ are connected by a Reidemeister-$i$ move means that (1) both of the diagrams are in  $S(c,d)$, and (2) one of $D,D'$ is changed to the other by a Reidemeister-$i$ move.

\noindent $O_{c,d}$ means that the value is invariant under changing order of components for any two marked link diagrams in $S(c,d)$.

\noindent $All_{c}$ means that the value is invariant under changing of base points, order of components and Reidemeister moves for any oriented link diagrams with at most $c$ crossings.

\noindent $D_{c}$ means that the value for any marked link diagram with $\leq c$ crossings is uniquely defined.

\noindent Similarly, we have $R^i_c, B_{c}$, $O_{c}$, and etc.

\end{remark}

\begin{proof}

The construction and proofs are all using induction on the index pair $(c,d)$, where $c$ is the crossing number of the diagram, and $d$ is the number of bad points of the diagram. It is obvious that $0\leq d\leq c$.

$\\$
\noindent {\bf The initial Step.} For a diagram of index $(c,0)$, namely a monotone diagram with $c$ crossing points, define its value to be $v_n$, where $n$ mean that the link has $n$ components.

\noindent Then the statement (0)-(4) is satisfied for diagrams in $S(c,0)$. This means that, for example, if $D,D'\in S(c,0)$ and the difference between them is a base point change, then $f(D)=f(D')$. Hence we have $All_0$.

$\\$
\noindent {\bf The inductive Step.}
Now suppose the statement (0)-(4) is proved for link diagrams with indices strictly less than $(c,d)$, in particular, we have $All_{c-1}$. This mean for any marked oriented link diagram with crossings $<c$, the value of the invariant is uniquely defined, independent of choice of base points and ordering of link components. Hence we can choose base points and ordering of link components arbitrarily for it to define the invariant.

$\\$
\noindent {\bf Proof of the statement (0)}: $\{All_{c-1},D_{c,d-1}\}\Rightarrow D_{c,d} \Rightarrow Sf_{c,d}$.

If the diagram $D$ has bad points, say its index is $(c,d)$, where $d>0$, we resolve the diagram at the first bad point $p$. Then, in the corresponding skein equation, all the other terms are of smaller indices than $(c,d)$.

For those diagrams with smaller crossing number, their orientations are given by the skein relation, we need to arbitrarily choose base points and ordering of link components. Then the invariant is uniquely defined, and by induction hypothesis $All_{c-1}$, the definition is independent of choice of base points and ordering of link components. Hence their values are uniquely defined. There is one term corresponds to crossing change, and it has a canonical orientation, base point set and ordering of link components (same as $D$). It also has one less bad points, hence it is already defined by induction ($D_{c,d-1}$). So all the terms except $D$ in the skein equation have been uniquely defined. Now we ask $b,b'$ both have left inverses. Hence the skein relation uniquely determines the value for $D$. We take this as the definition of invariant for $D$. Hence we have $Sf_{c,d}$. We shall prove later that if we resolve at other crossing point we shall get the same result.

\begin{remark}
We can similarly define the invariant for marked diagrams on $S^2$. Given a marked link diagram $D$ on $R^2$, we can also regard it as a marked diagram on $S^2$. However, for marked link diagram $D$ on $S^2$, we can have many marked diagram on $R^2$, depending on where we pick the $\infty$ point. All those marked diagrams on $R^2$ have the same value of invariant by the definition above. As a consequence, when we later prove the Reidemeister moves invariance, we can actually allow more "generalized Reidemeister moves". For example, if an outermost monogon contains the $\infty$ point, we can use the Reidemeister move I to reduce it.
\end{remark}

\noindent {\bf Proof of the statement (1)}: $\{All_{c-1},D_{c}\}\Rightarrow Sb_{c,d}$.

For a link diagram $D$, if $D$ has one bad point, then by definition, it satisfies the statement (1).
If $D$ has at least 2 bad points, and one resolve at a bad point $q$. If $q$ is the first bad point, then by definition, the equation is satisfied.
If not, denote the first bad point by $p$, and denote the value of $D$ by $f(D)$. If we resolve at $p$, we get many diagrams $D_1,D_2,\cdots $ and a linear sum $f_p(D)=\sum \alpha_i f(D_i)$ for some $\alpha_i$. Then by definition $f(D)=f_p(D)$.

We resolve each $D_i$ at $q$, then we get the linear sum $f_{q}(D_i)$. Each diagram $D_i$ has strictly lower indices than $(c,b)$. If $D_i$ has crossing number $c-1$, then skein equation is proved for resolving at any point. If $D_i$ has crossing number $c$, then it has $b-1$ bad points, and $q$ is also a bad point of $D_i$. In both cases, by induction hypothesis, $f(D_i)=f_{q}(D_i)$. Hence $f(D)=f_p(D)=\sum \alpha_i f_{q}(D_i)$.

On the other hand, we can resolve $D$ at $q$ first, we get many diagrams $D_1',D_2',\cdots $, each has strictly lower indices than $(c,b)$. Hence the statements (0)-(4) are satisfied. We get a linear sum $f_q(D)=\sum \beta_i f(D_i')$. We resolve each $D_i'$ at $p$, then we get the linear sum  $f_p(D_i')$. By the argument before and our induction hypothesis, $f(D_i')=f_p(D_i')$. Hence $f_q(D)=\sum \beta_i f_p(D_i')$. On the other hand, the ring is designed such that $\sum \beta_i f_p(D_i')=\sum \alpha_i f_{q}(D_i)$! (This is the equation $f_{pq}=f_{qp}$.)
$$\xymatrix{f(D) \ar@{=} [d]^{definition}_{1st\ bad\ point} &    \\ f_p(D)=\sum \alpha_i f(D_i) \ar@{=} [d]_{induction}^{hypothesis} & f_q(D)=\sum \beta_i f(D_i') \ar@{=} [d]_{induction}^{hypothesis} \\ \sum \alpha_i f_{q}(D_i) \ar@{=} [r]^{f_{pq}=f_{qp}} & \sum \beta_i f_p(D_i') }$$

Therefor, $f(D)=f_p(D)=\sum \alpha_i f_{q}(D_i)=\sum \beta_i f_p(D_i')=f_q(D)$. That is, if we resolve at $q$, the skein equation is satisfied.

\begin{cor} If one resolve at any point (not necessarily bad), the  skein equation is satisfied. $\{All_{c-1},D_{c}\}\Rightarrow S_c$.
\end{cor}

\begin{proof} If $q$ is a good point of $D$, we make a crossing change at $q$ get a new diagram $D'$, then $q$ is bad point of $D'$. The above proves that if we resolve $D'$ at $q$ the skein equation is satisfied. But this the same equation of $D$ resolving at $q$.

\end{proof}

\noindent This means that one can resolve at any crossing point to calculate the invariant, not necessarily the first bad point.

\noindent {\bf Proof of the statement (2) and (3)}:

\noindent This the hardest part and key of the whole proof.

\begin{lem} $\{All_{c-1},D_{c}, S_c \}\Rightarrow $ good $R^1_c$ and good $R^2_c$.
\end{lem}

\noindent Here, good $R^1_c$ means that for any two marked diagrams $D,D'$ with $\leq c$ crossings, if they are connected by a Reidemeister one move, and the crossing point $p$ in the Reidemeister one move is a good point, then $f(D)=f(D')$. Similarly, we can define the good $R^2_c$.

\begin{proof} The two diagrams $D,D'$ have the same number of bad points, one to one corresponds to each other. If there exists a bad point $q$ in both $D,D'$, we use the the skein relation, i.e. $S_c$, to resolve the bad point, then we have a pair of diagrams $\overline{D},\overline{D}'$ corresponds to crossing change, and other with lower crossing number. The diagrams with lower crossing number are equal to each other in pairs by induction on crossing number, since we can change the base point so that $p$ is still a good point. Hence $f(D)=f(D')$ if and only if $f(\overline{D})=f(\overline{D}')$. Since $\overline{D},\overline{D}'$ have fewer bad points than $D,D'$, this reduces to the case where $D,D'$ have no bad points at all. In this case, we have two monotone diagrams, hence by definition $f(D)=f(D')=v_n$, where $n$ denotes the number of link components of $D$.

\noindent Similarly, we can prove the good $R^2_c$.

\end{proof}

\begin{remark}
In this proof, removing bad points using $S_c$ and them reduce by induction is a key technique for proving our results. We shall simply refer to it as "remove the other bad points".
\end{remark}

\noindent $\{All_{c-1},D_{c}, S_c \}\Rightarrow $  $R^3_c$.

\begin{proof} Given two diagrams $D$ and $D'$, which differs by a Reidemeister move III.  Like above, we can assume all other points are good. In the two local diagrams containing the Reidemeister move III, there is a one to one correspondence between the three arcs appearing in the two local diagrams. We can order the three arcs by 1,2,3,($1',2,3'$ in $D'$) such that arc 1 ($1'$) is above arc 2 ($2'$), and arc 2 ($2'$) is above arc 3 ($3'$). The one to one correspondence preserves the ordering. Their intersections induce a one to one correspondence between the three pair points in the two diagrams.  Call them $p,p'$, $q,q',r,r'$. If arc $i$ intersects arc $j$ at $x$, then arc $i'$ intersects arc $j'$ at $x'$.

Suppose $p$ is the intersection of arc 1 and arc 2 (or arc 2 and arc 3), then we can resolve both $p$ and $p'$. Then we get many new link diagrams. There is a canonical one to one correspondence between those diagram, so we can denote them by $D_1,D_2, \cdots , D_1',D_2', \cdots $. Here $D_1,D_1'$ correspond to crossing change for $D$ and $D'$, and all other diagrams are of smaller crossing numbers. For those diagrams, we have $f(D_i)=f(D_i')$, $i\geq 2$ (By $R^2_{c-1}$ and $All_{c-1}$). Therefor, $f(D)=f(D')$ if and only if $f(D_1)=f(D_1')$. So we can assume $p$ is a good point. Similarly, we can assume the intersection of arc 2 and arc 3 is a good point.

Now, the intersection of arc 1 and arc 3, say $r$, is also a good point. The reason is simple. Since we proved base point invariance, we can assume there is no base point on any of the 3 arcs. The intersection of arc 2 and arc 3 is good means we first travel arc 3, then arc 2. Likewise, intersection of arc 1 and arc 2 is good means we first travel arc 2, then arc 1. Hence we first travel arc 3, then arc 1. Hence the intersection of arc 1 and arc 3 is good.

Hence one can make all the three intersections $p,q,r$ good. It follows that $p',q',r'$ are good. Now we have two monotone diagrams, the invariance is clear.
\end{proof}

\begin{lem} (\cite{MR1472978} Lemma 15.1) Suppose that $p$ and $q$ are two arcs in $R^2$ meeting only at their end points $A$ and $B$, and let $R$ be the compact region bounded by $p\cup q$. Suppose that $t_1, t_2, \cdots , t_n$ are arcs in $R$, each meeting $p\cup q$ at just its end points, one in $p$ and one in $q$. Suppose that every $t_i \cap t_j$ is at most one point, that intersections of arcs are transverse and there are no triple points. The graph, with vertices all intersections of these arcs and edges comprising $p\cup q \cup (\cup_i t_i )$, separates $R$ into collection of $v$-gons. Then amongst these  $v$-gons there is a 3-gon with an edge in $p$ and a $v$-gons there is a 3-gon with an edge in $q$.
\end{lem}

Using the above lemma, and a modification of {\cite{Manturov2004}} Lemma 5.1, we can prove the following lemma for link diagrams on $S^2$.

\begin{lem}\label{lem:3} (1) Each marked link diagram $D$ on $S^2$ with $\leq c$ crossings can be transformed to the unlink diagram without crossing by the following operation: crossing change, good $R^1_{c-1}$, good $R^2_{c-1}$, and $R^3_c$.

\noindent (2) Furthermore, for any give crossing point $p$ in the diagram $D$, we can use the above operations to  transform $D$ to the unlink diagram with only crossing $p$ such that the operations do not involve $p$.

\noindent (3) Also, if $D$ can be reduce by Reidemeister II move, where the two crossings in the Reidemeister II move are $p, q$, then  we can use the above operations to transform $D$ to the unlink diagram with only two crossings $p,q$, such that the operations do not involve $p,q$.
\end{lem}

\begin{proof}
The proof is almost the same as in {\cite{Manturov2004}}, except that other than the innermost argument there, we can also use an outermost argument to remove an bigon or monogon that contains the $\infty$ point using good $R^1_{c-1}$ or good $R^2_{c-1}$.
\end{proof}

\noindent $\{All_{c-1},D_{c}, S_c \}\Rightarrow B_c$.

Given a diagram $D$ with a fixed orientation and order of components, suppose that there are two base point sets $B$ and $B'$. We only need to deal with the case that $B$ and $B'$ has only one point $x$ and $x'$ different, they are in the same component $k$, and between $x$ and $x'$ there is only one crossing point $p$. In the base point sets $B$ and $B'$, $D$ has the same bad points except $p$. We shall prove the equation $f_B(D)=f_{B'}(D)$. If there is bad point other than $p$, say $q$, we resolve $D$ at $q$ to get diagrams $D_1,D_2,\cdots $. Then those $D_i$'s has lower indices than $D$, hence base point invariance is proved for them. As before, we get a marked diagram $\overline{D}$ corresponding to crossing change at $q$, and  $f_B(D)=f_{B'}(D)$ if and only if  $f_B(\overline{D})=f_{B'}(\overline{D})$. Hence we can assume there are no other bad points.

If there is no other bad points, there are two cases. Case 1. $p$ is a good point for both the two base point systems, then the values for $D$ are both $v_n$, hence equal.
Case 2. $p$ is a bad point for both the two base point systems, then the skein equation tells the values are the same.

Case 3, $p$ is good in $B$, bad in $B'$.
Then the diagram $D$ with base point set $B$ is a monotone diagram. Applying the above lemma~\ref{lem:3}(1) to $D$, we can fix the crossing $p$ and reduce all other crossings by crossing change, good $R^1_{c-1}$, good $R^2_{c-1}$, and $R^3_c$. The result is a link diagram $D'$ with only one crossing, $p$, and $f_B(D)=f_{B'}(D)$ if and only if  $f_B(D')=f_{B'}(D')$.

It follows that for $D'$, all the smoothings at $p$ (using skein relation) produce trivial links.
Then, with base point set $B$, the value of $D'$ is $v_n$, since it is a monotone diagram. In $B'$, the value is uniquely defined by the skein equation. Suppose the value with base point set $B'$ is $v_n$, then plug this into the skein equation, we necessarily have $(1+b+d_1+d_2)v_n+(c_1+c_2+c_3+c_4)v_{n+1}=0$.

On the other hand, this is also a sufficient condition. This because that the skein equation defines the value of the diagram. So, as long as the symbols always satisfy the equation $(1+b+d_1+d_2)v_n+(c_1+c_2+c_3+c_4)v_{n+1}=0$ for any $n\geq 1$, the value for $D'$ with base point set $B'$ is $v_n$. This proves the base point invariance.

\noindent $\{All_{c-1},D_{c}, S_c, B_c \}\Rightarrow R^1_c$.

\noindent {\bf (i)} Given two diagrams $D$ and $D'$, which differ by a Reidemeister one move. Say $D$ has index $(c,d)$, where $D'$ has index $(c+1,d')$. $D'$ has an extra crossing point $p$ because of the Reidemeister one move. $D$ and $D'$ have the same bad points except $p$. Like in previous proofs, if there are bad points other than $p$, we can resolve them and prove Reidemeister move one invariance inductively.

Otherwise, all other points are good, then $D$ and $D'$ are both diagrams of trivial links. If $p$ is good in $D'$, there is nothing to prove. If
$p$ is bad, we can use $B_c$, base points invariance, to get rid of this bad point by changing base point and then get the proof.

\noindent $\{All_{c-1},D_{c}, S_c, B_c \}\Rightarrow R^2_c$.

\noindent {\bf (ii)} Given two diagrams $D$ and $D'$, which differs at a Reidemeister move II. $D'$ has two more crossings, $p$ and $q$. Likewise, we can assume all other points are good. If the two crossings, $p$ and $q$, one is good, the other is bad, one can use a base point change to make them both good. Then both the diagrams $D$ and $D'$ are monotone diagrams. There is nothing to prove.

The only case needs proof is that both the two crossing are bad, and base point changes wouldn't change them from bad to good. However, changing both the two crossing will make them both good. Hence both the diagrams are diagrams for trivial link. In this case, we can apply lemma~\ref{lem:3}(3) to reduce all other crossings in the diagram. Hence we have the case that one diagram $D$ is a trivial (this means there are no crossing points), the other $D'$ has only two bad crossings. The crossings are intersections from two link components (otherwise we can use Reidemeister move I to reduce it). We have the following Fig. ~\ref{f6}.
\\
\\
\\
\\
\\
\\
\\
\\
\\

\begin{figure}[h]
\begin{center}
\scalebox{0.6}
{
\psset{arrowscale=2,unit=0.6}
\rput[0,0](0,9)
{
\rput[0,0](0,0)
{
\begin{pspicture}
(0,0)(18,0)
\psarc[linewidth=1pt]{->}(0,1.5){2}{-70.0}{70.0}
\psarc[linewidth=1pt]{<-}(3,1.5){2}{110.0}{130.0}
\psarc[linewidth=1pt](3,1.5){2}{150.0}{215.0}
\psarc[linewidth=1pt](3,1.5){2}{225.0}{250.0}
\usefont{T1}{ptm}{m}{n}
\rput(3,1.5){$X_1$}
\rput(2,2.7){$-$}
\rput(2,0.3){$+$}
\end{pspicture}
}
\rput[0,0](5,0)
{
\begin{pspicture}
(0,0)(18,0)
\psarc[linewidth=1pt]{->}(0,1.5){2}{-70.0}{70.0}
\psarc[linewidth=1pt](3,1.5){2}{110.0}{130.0}
\psarc[linewidth=1pt](3,1.5){2}{150.0}{215.0}
\psarc[linewidth=1pt]{->}(3,1.5){2}{225.0}{250.0}
\usefont{T1}{ptm}{m}{n}
\rput(3,1.5){$X_2$}
\rput(2,2.7){$+$}
\rput(2,0.3){$-$}
\end{pspicture}
}
\rput[0,0](10,0)
{
\begin{pspicture}
(0,0)(18,0)
\psarc[linewidth=1pt]{<-}(0,1.5){2}{-70.0}{70.0}
\psarc[linewidth=1pt]{<-}(3,1.5){2}{110.0}{130.0}
\psarc[linewidth=1pt](3,1.5){2}{150.0}{215.0}
\psarc[linewidth=1pt](3,1.5){2}{225.0}{250.0}
\usefont{T1}{ptm}{m}{n}
\rput(3,1.5){$X_3$}
\rput(2,2.7){$+$}
\rput(2,0.3){$-$}
\end{pspicture}
}
}
\rput[0,0](0,4)
{\rput[0,0](0,0)
{
\begin{pspicture}
(0,0)(18,0)
\psarc[linewidth=1pt](0,1.5){2}{-70.0}{-55.0}
\psarc[linewidth=1pt](0,1.5){2}{-30.0}{30.0}
\psarc[linewidth=1pt]{->}(0,1.5){2}{50.0}{70.0}
\psarc[linewidth=1pt]{<-}(3,1.5){2}{110.0}{250.0}
\usefont{T1}{ptm}{m}{n}
\rput(3,1.5){$X_1'$}
\rput(2,2.7){$+$}
\rput(2,0.3){$-$}
\end{pspicture}
}
\rput[0,0](5,0)
{
\begin{pspicture}
(0,0)(18,0)
\psarc[linewidth=1pt](0,1.5){2}{-70.0}{-55.0}
\psarc[linewidth=1pt](0,1.5){2}{-30.0}{30.0}
\psarc[linewidth=1pt]{->}(0,1.5){2}{50.0}{70.0}
\psarc[linewidth=1pt]{->}(3,1.5){2}{110.0}{250.0}
\usefont{T1}{ptm}{m}{n}
\rput(3,1.5){$X_2'$}
\rput(2,2.7){$-$}
\rput(2,0.3){$+$}
\end{pspicture}
}
\rput[0,0](10,0)
{
\begin{pspicture}
(0,0)(18,0)
\psarc[linewidth=1pt]{<-}(0,1.5){2}{-70.0}{-50.0}
\psarc[linewidth=1pt](0,1.5){2}{-30.0}{30.0}
\psarc[linewidth=1pt](0,1.5){2}{50.0}{70.0}
\psarc[linewidth=1pt]{<-}(3,1.5){2}{110.0}{250.0}
\usefont{T1}{ptm}{m}{n}
\rput(3,1.5){$X_3'$}
\rput(2,2.7){$-$}
\rput(2,0.3){$+$}
\end{pspicture}
}}
\rput[0,0](0,-1)
{\rput[0,0](0,0)
{
\begin{pspicture}
(0,0)(18,0)
\psarc[linewidth=1pt](0,1.5){2}{-70.0}{-55.0}
\psarc[linewidth=1pt]{->}(0,1.5){2}{-30.0}{70.0}
\psarc[linewidth=1pt]{<-}(3,1.5){2}{110.0}{130.0}
\psarc[linewidth=1pt](3,1.5){2}{150.0}{250.0}
\usefont{T1}{ptm}{m}{n}
\rput(3,1.5){$Y_1$}
\rput(2,2.7){$-$}
\rput(2,0.3){$-$}
\end{pspicture}
}
\rput[0,0](5,0)
{
\begin{pspicture}
(0,0)(18,0)
\psarc[linewidth=1pt](0,1.5){2}{-70.0}{35.0}
\psarc[linewidth=1pt]{->}(0,1.5){2}{50.0}{70.0}
\psarc[linewidth=1pt](3,1.5){2}{110.0}{210.0}
\psarc[linewidth=1pt]{->}(3,1.5){2}{230.0}{250.0}
\usefont{T1}{ptm}{m}{n}
\rput(3,1.5){$Y_2$}
\rput(2,2.7){$-$}
\rput(2,0.3){$-$}
\end{pspicture}
}
\rput[0,0](10,0)
{
\begin{pspicture}
(0,0)(18,0)
\psarc[linewidth=1pt]{<-}(0,1.5){2}{-70.0}{35.0}
\psarc[linewidth=1pt](0,1.5){2}{50.0}{70.0}
\psarc[linewidth=1pt]{<-}(3,1.5){2}{110.0}{210.0}
\psarc[linewidth=1pt](3,1.5){2}{230.0}{250.0}
\usefont{T1}{ptm}{m}{n}
\rput(3,1.5){$Y_3$}
\rput(2,2.7){$-$}
\rput(2,0.3){$-$}
\end{pspicture}
}}
}
\end{center}
\caption{Reidemeister II invariance.}\label{f6}
\end{figure}
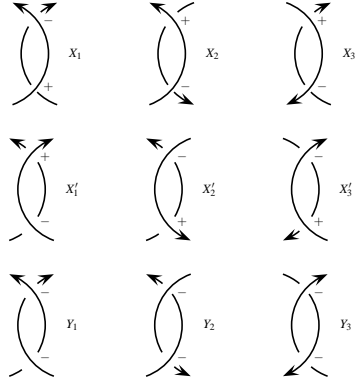

Those diagrams contain all cases of Reidemeister moves. Look at each column. We shall show that $f(X_i)=f(X_i')$ for $i=1,2,3$.
Let's resolve both $X_i$ and $X_i'$ at the positive crossing point, then we have

$X_i+bY_i+(c_1'+c_2')v_{n-1}+(d_1'+d_2')v_{n-1}=0$ and $X_i'+bY_i+(c_1'+c_2')v_{n-1}+(d_1'+d_2')v_{n-1}=0$. Hence  $f(X_i)=f(X_i')$ for $i=1,2,3$.

By changing from $X_i$ to $X_i'$ or from $X_i'$ to $X_i$, we can reduce the number of bad points by 2. Hence we get a monotone diagram. Hence we have $f(X_i)=f(X_i')=v_n=f(D)$.

\bigskip
\noindent {\bf Proof of the statement (4)}: \noindent $\{All_{c-1},D_{c}, S_c, B_c, R_c \}\Rightarrow O_c$.

\noindent Given two marked diagrams with different ordering of components. For simplicity, call them $D^1$ and $D^2$. By lemma~\ref{lem:3}(1), they can be simultaneously reduced to trivial marked diagrams $\overline{D}^1$ and $\overline{D}^2$ by crossing change, good $R^1_{c-1}$, good $R^2_{c-1}$ and $R^3_c$. We have $f(D_1)=f(D_2)$ if and only if $f(\overline{D}_1)=f(\overline{D}_2)$. However, $\overline{D}^1$ and $\overline{D}^2$ are trivial link diagrams with different ordering of link components. By definition, $f(\overline{D}_1)=f(\overline{D}_2)$. Hence $f(D_1)=f(D_2)$.

\end{proof}

\begin{cor}
The last set of relations $R^{A_1}_3$ for $A_1$ is $(1+b+d_1+d_2)v_n+(c_1+c_2+c_3+c_4)v_{n+1}=0$ for base point invariance.
\end{cor}

\subsection{Modifying by writhe}\label{sec:2.3}

There is another closely related link invariant with values in another commutative ring $A'_2$. The idea is that the skein relations can reduce the calculation to monotone diagrams, and we can regard the set of monotone diagrams as a basis and assign writhe dependant values to those diagrams. Now the skein relations don't give us a link invariant, but we can make a new function $g(w)$, such that the product $g(w)f(D)$ is regular link invariant. Here $w$ is a the writhe of the link diagram, regular means it is invariant under Reidemeister move tow and three. This construction is similar to the Kauffman bracket and the Kauffman 2-variable polynomial.

\begin{prop} There are functions $f,F$ for oriented link diagrams satisfying the following properties.

\noindent (1) For a monotone diagram $D$, $f(D)=h(w)v_n$, $F(D)=v_n$, where $w$ is a the writhe of the link diagram, $n$ is the number of components. There is another function $g(w)$ such that  $F(D)=g(w)f(D)$.

\noindent (2) For any marked link diagram $D$, $f(D)$ and $F(D)$ are uniquely defined.

\noindent (3) The function $f$ satisfies type one skein relations if we resolve at any bad point.

\noindent (4) The functions $f,F$ are invariant under base point change.

\noindent (5) $F$ and $f$ are invariant under Reidemeister moves II and III.

\noindent (6) $F$ and $f$ are invariant under changing order of components.

\end{prop}

\noindent We can regard $w$ as a function on oriented diagrams with integral values. Then we can write the skein relations for $F$ as follows.  For $f(E_{+})+bf(E_{-})+c_1f(E)+c_2f(W)+c_3f(HC)+c_4f(HT)+d_1f(VC)+d_2f(VT)=0$, since $f(D)=h(w(D))F(D)$, we can get $h(w(E_{+}))F(E_{+})+bh(w(E_{-}))F(E_{-})+c_1h(w(E))F(E)+c_2h(w(W))F(W)+c_3h(w(HC))F(HC)+c_4h(w(HT))F(HT)+d_1h(w(VC))F(VC)+d_2h(w(VT))F(VT)=0$.
Similarly, we have $h(w(E_{+}))F(E_{+})+b'h(w(E_{-}))F(E_{-})+c_1'h(w(E))F(E)+c_2'h(w(W))F(W)+d'_1h(w(S))F(S)+d'_2h(w(N))F(N)=0$.

\begin{remark}
As before, we shall prove it inductively and use the following notations.

\noindent $D^f_{c,d}$ means that the value of $f$ for any marked link diagram in $S(c,d)$ is uniquely defined, and $D^F_{c,d}$ means that the value of $F$ for any marked link diagram in $S(c,d)$ is uniquely defined. $D_{c,d}$ means that both $D^f_{c,d}$ and $D^F_{c,d}$ are true.  Similarly, we have $Sf^f_{c,d}$ and $Sf^F_{c,d}$, and etc.
\end{remark}

\begin{proof}
As before, the proof is an induction on index $(c,d)$. For statement (1), there is nothing to prove. We shall tell more about $g(w)$ later.

\bigskip\noindent {\bf Proof of the statement (2)(3)}: The proof is almost the same as in last section. We also resolve at the first bad point to define the invariant.

\bigskip\noindent {\bf Proof of the statement (4)(5)}:

We shall prove  by induction. Suppose that we have $All_{c-1}, D_c, S_c$. For diagrams in $S(c)=\cup_{d\leq c} S(c,d)$, we shall prove the following lemma first.

\begin{lemma}\label{lem:3}
Given two marked link diagrams $D,D'$ which are connected by a good Reidemeister one move, where $D'$ has one more crossing $p$ than $D$, then there is a shifting operator $T$, such that if $p$ is positive, then $f(D')=Tf(D)$, if $p$ is negative, then $f(D')=T^{-1}f(D)$.
\end{lemma}

\begin{proof}
We resolve $D,D'$ simultaneously at all bad points other than $p$, and for the resulting diagrams, we resolve again all bad points other than $p$, and so on. Finally, $f(D),f(D')$ can be written as linear combination of monotone diagrams $\{D_i,D_i'\mid i=1, \cdots ,N\}$. Those diagrams can be paired up such that $D_i,D_i'$ are connected by a good Reidemeister one move, where $D_i'$ has one more crossing $p$ than $D_i$. Suppose $p$ is positive. Then by definition, $f(D_i')=h(w(D_i)+1)v_n$, $f(D_i)=h(w(D_i))v_n$. Define $T(h(w))=h(w+1)$, then $f(D_i')=Tf(D_i)$ for each $i$. We can linearly extend the definition of $T$, hence $f(D')=Tf(D)$.

\noindent Likewise, if $p$ is negative, then $f(D')=T^{-1}f(D)$.
\end{proof}

\noindent $All_{c-1},D_c,S_c \Rightarrow B_c$.

\noindent Suppose we have a diagrams with different base point sets $B,B'$. For simplicity, the diagram will be called $D$ with base point set $B$, and $D'$ with base point set $B'$. $B$ and $B'$ has only one point $x$ and $x'$ different, they are in the same component $k$, and between $x$ and $x'$ there is only one crossing point $p$. In the base point systems $B'$, $D'$ has one bad point $p$, and $D$ is a monotone diagram.

Using the skein relations for $f(D)$, we can prove similar results as for type one invariants. For example, good $R^2_c$, $R^3_c$. We use lemma~\ref{lem:3} to replace good $R^1_c$. We also can simultaneously resolve any crossing points. This means that for $D,D'$, we can simultaneously apply crossing change, good $R^1_c$, good $R^2_c$, and $R^3_c$. Let the results be $\overline{D},\overline{D}'$. Then $f(D)=f(D')$ if and only if $f(\overline{D})=f(\overline{D}')$. The proof for type one invariant also applies here.

For example, if we have two diagrams $D,D'$ for both which we can apply a good Ridemeister one move to remove one good crossing $p$, the two oriented diagrams are the same but with different base point sets, then in the skein equation for any crossing point, all the lower crossing terms by induction ($D_{c-1}, B_{c-1}$ and $O_{c-1}$) are equal in pairs. Then $f(D)=f(D')$ if and only if $f(\overline{D})=f(\overline{D}')$. Hence we can make crossing changes at all crossings.

Hence the diagrams $D,D'$ are reduced to have only one crossing point $p$.

If $p$ is a positive crossing, then $f(\overline{D})=h(1)v_n$, $f(\overline{D}')+h(-1)bv_n+(c_1+c_2+c_3+c_4)h(0)v_{n+1}+(d_1+d_2)h(0)v_n=0$. Hence we need the equation $h(1)v_n+h(-1)bv_n+(c_1+c_2+c_3+c_4)h(0)v_{n+1}+(d_1+d_2)h(0)v_n=0$.

Similarly, if the crossing $p$ is negative, we get $f(\overline{D})=h(-1)v_n$, $h(1)v_n+bf(\overline{D}')+(c_1+c_2+c_3+c_4)h(0)v_{n+1}+(d_1+d_2)h(0)v_n=0$. Hence we need the equation $h(1)v_n+h(-1)bv_n+(c_1+c_2+c_3+c_4)h(0)v_{n+1}+(d_1+d_2)h(0)v_n=0$.

In both cases, we get $h(1)v_n+h(-1)bv_n+(c_1+c_2+c_3+c_4)h(0)v_{n+1}+(d_1+d_2)h(0)v_n=0$.

$f(D)=T^wf(\overline{D})$ for some integer $w$. Hence, for $f(D)=f(D')$ to be true, we need $h(w)v_n+h(w-2)bv_n+(c_1+c_2+c_3+c_4)h(w-1)v_{n+1}+(d_1+d_2)h(w-1)v_n=0$ hold in general.

Those equations guarantee base point sets invariance.

\noindent (ii) $All_{c-1},D_c,S_c, B_c \Rightarrow R^2_c$.

\noindent As before, we have two diagrams $D,D'$. We can assume that $D$ is a trivial link diagram, with $0$ crossing, and writhe is $0$. $D'$ has $2$ crossing. Those two bad crossings, say $p,q$ are intersections of different components. Then use the same argument as last proof for Reidemeister move II invariance, we get two equations. Last time, we had $X_i+bY_i+(c_1'+c_2')v_{n-1}+(d_1'+d_2')v_{n-1}=0$ and $X_i'+bY_i+(c_1'+c_2')v_{n-1}+(d_1'+d_2')v_{n-1}=0$. Now they should be modified a little bit, we have to add the writhe part into the equations. It is clear that the last 4 terms in the two equations all have writhe $-1$. So it is also true that $Xi$ and $Xi'$ have same value for $i=1,2,3$. After changing both $p,q$ to good points, the proof of Reidemeister move II invariance is trivial for $f$ and $F$.

\noindent (iii) $All_{c-1},D_c,S_c, B_c \Rightarrow R^3_c$.

Since Reidemeister move III does not change writhe, the proof is the same as before.

\bigskip\noindent {\bf Proof of the statement (6)} $All_{c-1},D_c,S_c, B_c, R_c \Rightarrow O_c$.

 Using lemma~\ref{lem:3}(1), the proof is the same as before.

\end{proof}

\begin{cor}
The last set of relations $R^{A_1'}_3$ for $A_1'$ contains the following relations:
$\\ h(w)v_n+h(w-2)bv_n+(c_1+c_2+c_3+c_4)h(w-1)v_{n+1}+(d_1+d_2)h(w-1)v_n=0.$
\end{cor}

The functions $f,F$ are not invariant under Reidemeister move one, with the help of base point sets invariance, we can discuss their behavior under Reidemeister move one. This time, we do not need to assume the  Reidemeister one move is a good one.

\begin{lemma}\label{lem:4}
Given two marked link diagrams $D,D'$ which are connected by a Reidemeister one move, where $D'$ has one more crossing $p$ than $D$, then there is a shifting operator $T$, such that if $p$ is positive, then $f(D')=Tf(D)$, if $p$ is negative, then $f(D')=T^{-1}f(D)$.
\end{lemma}

\begin{proof}
\noindent Given two marked link diagrams $D,D'$ which are connected by a Reidemeister one move, where $D'$ has one more crossing $p$ than $D$. By $B_c$, we can assume that $p$ is a good point. Now the result follows from lemma~\ref{lem:3}.
\end{proof}

Now, there are two ways to make link invariant out of this. The first, up to $T$ action equivalence, $f$ defines a link invariant. The second, we can use $F=g(w)f$ as follows. It is clear that $F$ is also a regular invariant, so we only need to worry about Reidemeister move one. Given two diagrams $D,D'$ connected by a Reidemeister one move, we have $F(D)=g(w)f(D)$ and $F(D')=g(w+1)Tf(D)$. To have $F(D)=F(D')$, we ask $g(w)f(D)=g(w+1)Tf(D)$. For this to be true, we interpret $g(w)$ as an operator on $A'_1$. $g(w)=T^{-w}$ will make the equation true! Hence $F(D)=T^{-w(D)}f(D)$ is a link invariant.

\begin{remark}
The difference with type one invariant is that although $F(D)=v_n$ on monotone diagrams, $F(D)$ does not satisfy the skein relations.
\end{remark}
\begin{remark}
\noindent An easy choice for the equation $h(w)v_n+h(w-2)bv_n+(c_1+c_2+c_3+c_4)h(w-1)v_{n+1}+(d_1+d_2)h(w-1)v_n=0$ is to let $h(w)=a^w$ for a new variable $a$. Then the equation is reduced to $a^2v_n+bv_n+(c_1+c_2+c_3+c_4)av_{n+1}+(d_1+d_2)av_n=0$. Then we can let $g(w)=a^{-w}$. From this, one can see that this is a new link invariant (If it is equivalent to type one invariant, we need to give some nontrivial proof).
\end{remark}

\section{The second link invariant}\label{sec:3}

   There is another closely related new link invariant.

   Given a link diagram with local crossing $E_{\pm}$, if the two strands are from same component, then they satisfies the following relation:

    $E_+=c_1E+c_2W+c_3HC+c_4HT+d_1VC+d_2VT$

    $E_-=\overline{c}_1E+\overline{c}_2W+\overline{c}_3HC+\overline{c}_4HT+\overline{d}_1VC+\overline{d}_2VT$

     \noindent Otherwise, they satisfies the following relation:

    $E_+=c'_1E+c'_2W+d'_1S+d'_2N$

    $E_-=\overline{c}'_1E+\overline{c}'_2W+\overline{d}'_1S+\overline{d}'_2N$

  Now, we have 20 variables. We can them the type two skein relations.
  Similarly, there is another set of equations for them, and the above setting defines a new link invariant.

\subsection{The ring $A_2$}\label{sec:3.1}

Since the skein relations of type one and type two invariants are so closely related, we can get the relation set for $A_2$ from the relation sets for $A_1$. It is very easy. Let $b=b'=0$ in relation set of $A_1$, then we get the relation set $R^{A_2}_2$ for $A_2$ as follows (they need to be complete by conjugation).

\bigskip\noindent {\bf Case 1:}  $c_2'\overline{c_2}'=\overline{d}'_1d_2'$, $d_1'c_2'+d_2'\overline{d}_2'=\overline{1}'c_2'+d_2'c_1'$.

\noindent {\bf Case 2:}  $c_1'c_4+c_2'c_4=c_3c_1'+c_3c_2'$,$c_1'd_1+c_2'd_1=\overline{d}_1'd_2'+\overline{d}_1'd_1'$, $c_2'd_2+c_1'd_2=\overline{d}_2'd_2'+\overline{d}_2'd_1'$, $d_2'\overline{c_4}+d_1'\overline{c_4}=\overline{c_3}d_1'+\overline{c_3}d_2'$.

\noindent {\bf Case 3:}  $c_4\overline{d}_1'=c_2c_2'$, $c_4\overline{2}'=c_2d_2'$, $c_2\overline{d}_2'+c_3c_2'=c_1c_2'+c_2c_1'$, $c_2\overline{c_1}'+c_3d_2'=c_1d_2'+c_2d_1'$, $c_2\overline{d}_1'=c_3c_2'$,  $c_2\overline{c_2}'=c_3d_2'$, $c_4\overline{d}_2'+c_1c_2'=c_4c_2'+c_4c_1'$, $c_4\overline{c_1}'+c_1d_2'=c_4d_2'+c_4d_1'$, $d_1\overline{d}_2'+d_2c_2'=d_1c_2'+d_1c_1'$, $d_1\overline{c_1}'+d_2d_2'=d_1d_2'+d_1d_2'$, $d_1\overline{d}_1'=d_2c_2'$, $d_1\overline{c_2}'=d_2d_2'$.

\noindent {\bf Case 4:}  $c_4\overline{d}_2'+c_3\overline{d}_1'=\overline{d}_2'c_3+\overline{d}_1'c_4$, $c_4\overline{d}_1'+c_3\overline{d}_2' =\overline{d}_1'c_3+\overline{d}_2'c_4$, $d_2\overline{d}_1+d_1\overline{d}_1=\overline{c_1}'c_3+\overline{c_2}'c_4+ d_2'c_2+d_1'c_1$, $d_1\overline{c_4}+d_2\overline{c_4}=\overline{c_1}d_2+\overline{c_2}d_1$,
$d_2\overline{d}_2+d_1\overline{d}_2=\overline{c_1}'c_4+\overline{c_2}'c_3+ d_1'c_2+d_2'c_1$, $d_1\overline{c_3}+d_2\overline{c_3}= \overline{c_1}d_1+\overline{c_2}d_2$, $c_3\overline{c_1}'+c_4\overline{c_2}'+ c_2d_2'+c_1d_2'=\overline{d}_1d_2+\overline{d}_1d_1$, $d_1\overline{c_1}+d_2\overline{c_2}=\overline{c_3}d_1+\overline{c_3}d_2$, $d_2\overline{c_1}+d_1\overline{c_2}=\overline{c_4}d_1+\overline{c_4}d_2$.

\noindent {\bf Case 5:}  $c_3c_2=c_4c_2$, $c_3c_3=c_2c_2$, $c_3d_2=c_2d_1$, $c_4c_2=c_2c_3$, $c_2c_4+c_4c_3=c_1c_3+c_2c_1$, $c_2d_1+c_4d_1=c_1d_1+c_2d_2$, $c_3c_1+c_1c_2=c_3c_2+c_3c_3$, $c_1c_4=c_3c_1$, $c_1d_1=c_4d_2$, $c_4c_3=c_2c_2$, $c_2c_3=c_4c_2$, $c_2d_2=c_4d_1$, $d_2c_2=d_1c_3$, $d_2c_3=d_1c_2$, $d_2d_2=d_1d_1$, $d_2c_1+d_1c_2=d_2c_2+d_2c_3$, $d_1c_4+d_2c_4=d_2c_3+d_2c_1$.

\subsection{Proofs for the second link invariant}\label{sec:3.2}

\begin{prop} The invariant satisfies the following properties.

(0) The value is defined uniquely for any oriented link diagram.

(1) Satisfying skein relations if we resolve at any crossing point.

(2) Invariant under Reidemeister moves for any two diagrams with crossing $<c$.

\end{prop}

$\\$
\noindent {\bf Proof of the statement (0)}:

We shall define the invariant inductively on crossing number $c$ of the diagram.

\noindent {\bf Step 1}. For a n component oriented link diagram of crossing number $c=0$, define its value to be $v_n$.

\noindent Then the statement (0)-(2) is satisfied.

\noindent {\bf step 2}. If the diagram $D$ has crossing points, we resolve the diagram at one crossing point $p$. Then, in the skein equation, all the other terms have smaller crossing numbers. By induction hypothesis, those terms are uniquely defined. Hence the skein relation uniquely defines the value for $D$.

\bigskip\noindent {\bf Proof of the statement (1)}:

The proof is almost the same as type 1 invariant.

\bigskip\noindent {\bf Proof of the statement (2)}:

\noindent {\bf (i)} Given two diagrams $D$ and $D'$, which differs at a Reidemeister move I. Say $D$ has index $c$, where $D'$ has index $c+1$. In the local Reidemeister move I part, say the crossing point is $p$. $D$ and $D'$ have the same crossing points except $p$. If there are crossing points other than $p$, we can resolve both the diagrams and prove Reidemeister move I invariance inductively.

Otherwise, $p$ is the only crossing point, then $D$ and $D'$ are both diagrams of trivial links. $f(D)$ is defined to be $v_n$, but $f(D')$ is defined by using the skein relation. After resolving at $p$, we get $f(D')=(c_1+c_2+c_3+c_4)v_{n+1}+(d_1+d_2)v_n$ or $f(D')=(\overline{c}_1+\overline{c}_2+\overline{c}_3+\overline{c}_4)v_{n+1}+(\overline{d}_1+\overline{d}_2)v_n$, depending on the crossing type.

Then\ Reidemeister move I invariance is guaranteed by the following equations
$v_n=(c_1+c_2+c_3+c_4)v_{n+1}+(d_1+d_2)v_n$ and $\\ v_n=(\overline{c}_1+\overline{c}_2+\overline{c}_3+\overline{c}_4)v_{n+1}+(\overline{d}_1+\overline{d}_2)v_n.$ We ask those equations to be always true.

\bigskip\noindent {\bf (ii)} Given two diagrams $D$ and $D'$, which differs at a Reidemeister II move. Likewise, we can assume there is no other crossing points. If the Reidemeister move II involves only one component, then it can be resolved by Reidemeister move I moves, and then invariance  followed from Reidemeister move I invariance.

If the Reidemeister move II involve two link components, we resolve it at the positive crossing point. Then Reidemeister move II invariance followed from the following equation: $ v_{n+1}=(\overline{c}_1+\overline{c}_2+\overline{d}_1+\overline{d}_2)v_{n}$.

\begin{remark}
If we resolve at the negative crossing point, we can get another invariant with a different relation   $v_{n+1}=(\overline{c}_1{}'+\overline{c}_2{}'+\overline{d}_1{}'+\overline{d}_2{}')v_{n}$.
\end{remark}

\bigskip\noindent {\bf (iii)} Given two diagrams $D$ and $D'$, which differs at a Reidemeister move III.  Like before, we can assume all other points are good. In the local diagram containing the Reidemeister move III, there is a one to one correspondence between the three arcs appearing in the two local diagrams. We can also order the three arcs by 1,2,3,($1',2,3'$ in $D'$) such that arc 1 is above arc, and arc 2 is above arc 3. The one to one correspondence preserve the ordering. Suppose arc 1 and arc 2 intersects at $p$, arc $1'$ and arc $2'$ intersects at $p'$. Then we can resolve at $p,p'$ at the same time. The resulting terms can be paired up and equal each other since we have proved Reidemeister move II invariance. Therefor, we proved Reidemeister move III invariance.

\subsection{Modifying by writhe}\label{sec:3.3}

Instead of asking the above definition to be Reidemeister moves invariant, we can ask its modification to be Reidemeister moves invariant. Denote the value of diagram $D$ by $f(D)$, let $\omega$ denote the writhe, $c$ denote the crossing number, $\mu$ denote number of link components. We ask a new family of functions $g(\omega , c, \mu)$ with parameters in  $\omega , c, \mu$, such that $gf(D)$ is Reidemeister moves invariant.

\begin{prop} There is an invariant satisfies the following properties.

\noindent (0) The value $f(D)$ is defined uniquely for any oriented link diagram.

\noindent (1) $f(D)$  satisfies skein relations if we resolve at any crossing point.

\noindent (2) $F(D)=g(\omega , c, \mu)f(D)$ is invariant under Reidemeister moves for any two diagrams with crossing $<c$.

\end{prop}

\begin{proof} We use inductions on crossing number $c$.

\bigskip\noindent {\bf Proof of the statement (0)(1)} Same as before.

\bigskip\noindent {\bf Proof of the statement (2)}

\noindent {\bf (i)} To make it Reidemeister move I invariant, like above, we need some new equations. Suppose there are two diagrams $D$ and $D'$, which differs at a Reidemeister move I. Say $D$ has index $c$, where $D'$ has index $c+1$. In the local Reidemeister move I part, say the crossing point is $p$. $D$ and $D'$ have the same crossing points except $p$. Suppose diagram $D$ has parameters $\omega , c, \mu$, the $D'$ has parameters $\omega +1 , c+1, \mu$ or has parameters $\omega -1, c+1, \mu$.

An easy case: $c=0$. If we use the skein relation to resolve the only crossing of $D'$, depending on the crossing type, we get
$D'=(c_1+c_2+c_3+c_4)v_{n+1}+(d_1+d_2)v_n$ or $D'=(\overline{c}_1+\overline{c}_2+\overline{c}_3+\overline{c}_4)v_{n+1}+(\overline{d}_1+\overline{d}_2)v_n$

In general, if they have other crossings, we can resolve all other crossings, and the resulting terms for $D$ and $D'$ can be paired up. Then we can resolve $p$. For example, if the $p$ has positive crossing, then we can group the terms for $D'$ together, such that each group has the form $(c_1+c_2+c_3+c_4)v_{n+1}+(d_1+d_2)v_n$ for some $n$, and for $D$, there is one term $v_n$ corresponds to it.

Hence we can add the following equations for Reidemeister move I invariance:

$g(\omega +1 , c+1, \mu)\{(c_1+c_2+c_3+c_4)v_{n+1}+(d_1+d_2)v_n\}=g(\omega , c, \mu)v_n$

$g(\omega -1, c+1, \mu)\{(\overline{c}_1+\overline{c}_2+\overline{c}_3+\overline{c}_4)v_{n+1}+(\overline{d}_1+\overline{d}_2)v_n\}=g(\omega , c, \mu)v_n$

\bigskip\noindent {\bf (ii)} To make it Reidemeister move II invariant, like above, we need some new equations. Suppose there are two diagrams $D$ and $D'$, which differs at a Reidemeister move II. Say $D$ has index $c$, where $D'$ has index $c+2$. In the local Reidemeister move II part, say the crossing point is $p,q$. $D$ and $D'$ have the same crossing points except $p,q$. Suppose diagram $D$ has parameters $\omega , c, \mu$, the $D'$ has parameters $\omega , c+2, \mu$.

We can resolve all other crossings, so that $f(D)$ and $f(D')$ turns to be linear combinations of diagrams $D_1,D_i',\cdots $ with coefficients. The resulting terms for $D$ and $D'$ can be paired up. Say $D_i$ and $D_i'$ is one of the pairs. $D_i$ has no crossings, and $D_i'$ has 2 crossings. We resolve $D_i'$ at the negative crossing point first, then the positive point, we always get

\noindent \resizebox{12.5cm}{!} {$(\overline{c}_1'+ \overline{c}_2')\{(c_1+c_2+c_3+c_4)v_{n}+(d_1+d_2)v_{n-1} \} +(\overline{d}_1'+ \overline{d}_2')\{(\overline{c}_1+\overline{c}_2+\overline{c}_3+\overline{c}_4)v_{n} +(\overline{d}_1+\overline{d}_2)v_{n-1}\}$}.

Then $f(D')$ is a linear combination, such that each term has the form $$(\overline{c}_1'+ \overline{c}_2')\{(c_1+c_2+c_3+c_4)v_{n}+(d_1+d_2)v_{n-1} \}+(\overline{d}_1'+ \overline{d}_2')\{(\overline{c}_1+\overline{c}_2+\overline{c}_3+\overline{c}_4)v_{n}+(\overline{d}_1+ \overline{d}_2)v_{n-1}\}$$
 for some $n$, and for $D$, there is one term $v_n$ corresponds to it.

Hence we can add the following equations for Reidemeister move II invariance:

\noindent $g(\omega , c+2, \mu)\{(\overline{c}_1'+ \overline{c}_2')\{(c_1+c_2+c_3+c_4)v_{n}+(d_1+d_2)v_{n-1} \}+(\overline{d}_1'+ \overline{d}_2')\{(\overline{c}_1+\overline{c}_2+\overline{c}_3+\overline{c}_4)v_{n}+ (\overline{d}_1+\overline{d}_2)v_{n-1}\} \}=g(\omega , c, \mu)v_n$

\bigskip\noindent {\bf (iii)} For the Reidemeister move III invariance, things are much easier. The only difference of $D$ and $D'$ is position of a crossing point. One choose that right point to resolve both the diagrams. The equality is trivial.

\end{proof}

\begin{cor}
$f(D)$ here is invariant under Reidemeister move III.

\end{cor}

\begin{remark} This is analogous to Kauffman's bracket and Jones polynomial.
\end{remark}

\section{Type 3 and type 4 invariants}\label{sec:4}

We can similarly define two more invariants.

\bigskip\noindent {\bf Type 3 invariant}

   Given a link diagram with local crossing $E_{\pm}$, if the two strands are from same component, then they satisfies the following relation:

   $E_{+}+bE_{-}+c_1E+c_2W+c_3HC+c_4HT+d_1VC+d_2VT=0$

     \noindent Otherwise, they satisfies the following relation:

    $E_+=c'_1E+c'_2W+d'_1S+d'_2N$

    $E_-=\overline{c}'_1E+\overline{c}'_2W+\overline{d}'_1S+\overline{d}'_2N$

\bigskip\noindent {\bf Type 4 invariant}

   Given a link diagram with local crossing $E_{\pm}$, if the two strands are from same component, then they satisfies the following relation:

    $E_+=c_1E+c_2W+c_3HC+c_4HT+d_1VC+d_2VT$

    $E_-=\overline{c}_1E+\overline{c}_2W+\overline{c}_3HC+\overline{c}_4HT+\overline{d}_1VC+\overline{d}_2VT$

     \noindent Otherwise, they satisfies the following relation:

    $E_{+}+b'E_{-}+c_1'E+c_2'W+d'_1S+d'_2N=0$

It seems likely that the above constructions also define some link invariants, but unfortunately, the equations $f_{pq}=f_{qp}$ tell us that in those equations many coefficients must be zero. The results are not of interest. In the future in our later paper, we will discuss some more general constructions, in which the ring is non commutative. We hope type 3 and 4 invariants will be non trivial in that case.

\section{Conclusion and theorems}\label{sec:5}

We have four rings here $A_1,A_1',A_2,A_2'$. We list below their generators and relation sets $R^{A_i}_1,R^{A_i}_3$ or $R^{A_i'}_1,R^{A_i'}_3$. Since $R^{A_i}_2=R^{A_i'}_2$ and they are too complicated, we do not list them here again.

$\\$
\noindent $A_1$ has generators $b,c_1,c_2,c_3,c_4,d_1,d_2,b',c_1',c_2',d_1',d_2'$ and $\{v_n\}_{n=1}^{\infty}$.

$R^{A_1}_1$: all generators commute.

$R^{A_1}_3$: $(1+b+d_1+d_2)v_n+(c_1+c_2+c_3+c_4)v_{n+1}=0$.

$\\$
\noindent $A_1'$ has generators $b,c_1,c_2,c_3,c_4,d_1,d_2,b',c_1',c_2',d_1',d_2'$ and $\{v_n\}_{n=1}^{\infty}$.

$R^{A_1'}_1$: all generators commute.

$R^{A_1'}_3$:
$h(w)v_n+h(w-2)bv_n+(c_1+c_2+c_3+c_4)h(w-1)v_{n+1}+(d_1+d_2)h(w-1)v_n=0$.

\begin{theorem}\label{thm:main} For oriented link diagrams, there is a link invariant with values in $A_1$ and satisfies the following skein relations:

\noindent (1) If the two strands are from same link component, then

$E_{+}+bE_{-}+c_1E+c_2W+c_3HC+c_4HT+d_1VC+d_2VT=0.$

\noindent (2) Otherwise,
$E_{+}+b'E_{-}+c_1'E+c_2'W+d'_1S+d'_2N=0.$

The value for trivial n-component link is $v_n$.

In general, replacing $A_1$ by any homomorphic image of $A_1$, one will get a link invariant.

There is a modified invariant taking values in $A_1'$, and the value for a monotone n-component link diagram is $h(w)v_n$.
\end{theorem}

This theorem directly follows from proposition 2.13. In general, if one replace $A_1$ by any homomorphic image of $A_1$, the equation $f_{pq}=f_{qp}$ still holds, then there is nothing to prove.

$\\$
\noindent $A_2$ has generators $c_1,c_2,c_3,c_4,d_1,d_2,c_1',c_2',d_1',d_2',\overline{c}_1,\overline{c}_2,\overline{c}_3,\overline{c}_4,\overline{d}_1, \overline{d}_2,\overline{c}_1',\overline{c}_2',\overline{d}_1',\overline{d}_2'$ and $\{v_n\}_{n=1}^{\infty}$.

\noindent $R^{A_2}_1$: all generators commute with each other.

\noindent $R^{A_2}_3$: $v_n=(c_1+c_2+c_3+c_4)v_{n+1}+(d_1+d_2)v_n$,

\noindent $v_n=(\overline{c}_1+\overline{c}_2+\overline{c}_3+\overline{c}_4)v_{n+1}+(\overline{d}_1+\overline{d}_2)v_n$,

\noindent $v_{n+1}=(\overline{c}_1+\overline{c}_2+\overline{d}_1+\overline{d}_2)v_{n}$.

$\\$
\noindent $A_2'$ has generators $c_1,c_2,c_3,c_4,d_1,d_2,c_1',c_2',d_1',d_2',\overline{c}_1,\overline{c}_2,\overline{c}_3,\overline{c}_4,\overline{d}_1, \overline{d}_2,\overline{c}_1',\overline{c}_2',\overline{d}_1',\overline{d}_2'$ and $\{v_n\}_{n=1}^{\infty}$.

\noindent $R^{A_2'}_1$: all generators commute with each other.

\noindent $R^{A_2'}_3$: $g(\omega +1 , c+1, \mu)\{(c_1+c_2+c_3+c_4)v_{n+1}+(d_1+d_2)v_n\}=g(\omega , c, \mu)v_n$,

\noindent $g(\omega -1, c+1, \mu)\{(\overline{c}_1+\overline{c}_2+\overline{c}_3+\overline{c}_4)v_{n+1}+(\overline{d}_1+\overline{d}_2)v_n\}=g(\omega , c, \mu)v_n$,

\noindent $g(\omega , c+2, \mu)\{(\overline{c}_1'+ \overline{c}_2')\{(c_1+c_2+c_3+c_4)v_{n}+(d_1+d_2)v_{n-1} \}+(\overline{d}_1'+ \overline{d}_2')\{(\overline{c}_1+\overline{c}_2+\overline{c}_3+\overline{c}_4)v_{n}+ (\overline{d}_1+\overline{d}_2)v_{n-1}\} \}=g(\omega , c, \mu)v_n$.

\begin{theorem}\label{thm:main2} For oriented link diagrams, there is a link invariant with values in $A_2$ and satisfies the following skein relations:

\noindent (1) If the two strands are from same link component, then

    $E_+=c_1E+c_2W+c_3HC+c_4HT+d_1VC+d_2VT$

    $E_-=\overline{c}_1E+\overline{c}_2W+\overline{c}_3HC+\overline{c}_4HT+\overline{d}_1VC+\overline{d}_2VT$

\noindent (2) Otherwise,

    $E_+=c'_1E+c'_2W+d'_1S+d'_2N$

    $E_-=\overline{c}'_1E+\overline{c}'_2W+\overline{d}'_1S+\overline{d}'_2N$

The value for a trivial n-component link diagram is $v_n$.

In general, replacing $A_2$ by any homomorphic image of $A_2$, one will get a link invariant.

There is a modified invariant taking values in $A_2'$, and the value for a trivial n-component link diagram is $v_n$.

\end{theorem}

\section{Relation with link polynomials}\label{sec:6}

Let's focus on type one invariant first. From the proof one can conclude that as long as any ring satisfies those relations, it will produce a link invariant. For example, if we let $c_2=c_3=c_4=d_1=d_2=0=c_2'=d_1'=d_2'$, and $b=b'$ then the relations are all gone. We get $E_++bE_-+c_1E=0$, this is equivalent to the well know HOMFLYPT polynomial! This means HOMFLYPT polynomial is a special case of the new invariants. By adding other relations into this relation set, we can design a special invariant. This is like group representations. You can add new relations to get a homomorphism image of the original group. We can refer to this as the homomorphisms of the invariants. In this section, we shall discuss some natural way to add new relations to get invariants which are easy to handle with. If we ask some symbols to be $0$,  and others to have inverses, it will produce many new link polynomials.

\bigskip\noindent (1) We can try the following: $c_1=c_2=c_3=c_4=0=c_1'=c_2'$, then the relations reduce to the following (need to be complete by conjugation)

\noindent {\bf Case 1:}  $0=\overline{d}'_1d_2'$, $b'd_2'=\overline{b}'d_2'$, $d_1'd_2'=d_2'd_1'$, $d_2'\overline{d}_2'=0$

\noindent {\bf Case 2:}  $d_2'\overline{b}=d_2'b'$, $d_1'\overline{b}=d_1'b'$, $0=\overline{d}_1'd_2'+\overline{d}_1'd_1'$, $0=\overline{d}_2'd_2'+\overline{d}_2'd_1'$.

\noindent {\bf Case 3:} $d_1\overline{d}_2'=0$, $d_2d_2'=d_1d_2'+d_1d_2'$, $d_1\overline{b}'=d_1b'$, $d_1\overline{d}_1'=0$, $0=d_2d_2'$.

\noindent {\bf Case 4:}  $d_2\overline{d}_1+d_1\overline{d}_1=0$,
$bd_2=\overline{b}d_2$, $d_2\overline{d}_2+d_1\overline{d}_2=0$, $bd_1=\overline{b}d_1$, $0=\overline{d}_1d_2+\overline{d}_1d_1$, $d_1\overline{b}=d_1b$, $d_2\overline{b}=d_2b$.

\noindent {\bf Case 5:} $d_2d_2=d_1d_1$

Those can be simplified to: $(b-\overline{b})d_1=0$, $(b-\overline{b})d_2=0$, $(b'-\overline{b}')d_1=0$, $(b'-\overline{b}')d_2=0$, $(b'-\overline{b}')d_1'=0$, $(b'-\overline{b}')d_2'=0$, $d_1d_1'=0$, $d_1d_2'=0$, $d_2d_2'=0$, $d_1'd_2'=0$, $d_1d_1=d_2d_2$, $d_1d_1+d_1d_2=0$.

Certainly, we need to complete them by conjugation. $R^{A_1}_3$ tells us that $(1+b+d_1+d_2)v_n=0$. Hence this ring has zero divisors.

An easy case is to let $b=b'=-1$, then we get the followings. $E_+-E_-+d_1VC+d_2VT=0$, $E_+-E_-+d_1'S+d_2'N=0$, $d_1d_1'=0$, $d_1d_2'=0$, $d_2d_2'=0$, $d_1'd_2'=0$, $d_1d_1=d_2d_2=-d_1d_2$, and $R^{A_1}_3=\emptyset $.

\bigskip\noindent (2)
Another homomorphism is to let $b=b'=1$. Then $x\equiv \overline{x}$. Simplify those relations, we get: $3=4$. The relations for $A_1$ can be simplified as followings:

\noindent {\bf Case 1:}  $c_2'c_2'=d'_1d_2'$, $d_1'c_2'+d_2'd_2'=c_1'c_2'+d_2'c_1'$.

\noindent {\bf Case 2:}  $c_1'd_1+c_2'd_1=d_1'd_2'+d_1'd_1'$, $c_2'd_2+c_1'd_2=d_2'd_2'+d_2'd_1'$, $d_2'c_3+d_1'c_3=c_3d_1'+c_3d_2'$.

\noindent {\bf Case 3:}  $c_3d_1'=c_2c_2'$, $c_3c_2'=c_2d_2'$, $c_2d_2'+c_3c_2'=c_1c_2'+c_2c_1'$, $c_2c_1'+c_3d_2'=c_1d_2'+c_2d_1'$, $c_2b=c_2b'$, $c_2d_1'=c_3c_2'$,  $c_2c_2'=c_3d_2'$, $c_3d_2'+c_1c_2'=c_3c_2'+c_3c_1'$, $c_3c_1'+c_1d_2'=c_3d_2'+c_3d_1'$, $d_1d_2'+d_2c_2'=d_1c_2'+d_1c_1'$, $d_1c_1'+d_2d_2'=d_1d_2'+d_1d_2'$, $d_1d_1'=d_2c_2'$, $d_1c_2'=d_2d_2'$.

\noindent {\bf Case 4:}  $d_2d_1+d_1d_1=c_1'c_3+c_2'c_3+ d_2'c_2+d_1'c_1$, $d_1c_3+d_2c_3=c_1d_2+c_2d_1$, $d_2d_2+d_1d_2=c_1'c_3+c_2'c_3+ d_1'c_2+d_2'c_1$, $d_1c_3+d_2c_3= c_1d_1+c_2d_2$, $c_3c_1'+c_3c_2'+ c_2d_2'+c_1d_2'=d_1d_2+d_1d_1$, $d_1c_1+d_2c_2=c_3d_1+c_3d_2$, $d_2c_1+d_1c_2=c_3d_1+c_3d_2$.

\noindent {\bf Case 5:}   $c_3c_3=c_2c_2$, $c_3d_2=c_2d_1$, $c_2c_3+c_3c_3=c_1c_3+c_2c_1$, $c_2d_1+c_3d_1=c_1d_1+c_2d_2$, $c_3c_1+c_1c_2=c_3c_2+c_3c_3$,$c_1d_1=c_3d_2$, $c_2d_2=c_3d_1$, $d_2c_2=d_1c_3$, $d_2c_3=d_1c_2$, $d_2d_2=d_1d_1$, $d_2c_1+d_1c_2=d_2c_2+d_2c_3$, $d_1c_3+d_2c_3=d_2c_3+d_2c_1$.

$R^{A_1}_3$ tells us that $(2+d_1+d_2)v_n+(c_1+c_2+c_3+c_4)v_{n+1}=0$. We can also let $b=b'=-1$, and get similar relations.

\bigskip\noindent (3) An interesting choice is to make the invariant like an "oriented version of Kauffman 2-variable polynomial" as follows. Let $c_2=c_3=c_4=d_2=c_2'=d_2'=0$, and $b=b'=1$, then we have the followings.

\noindent {\bf Case 2:}  $c_1'd_1=d_1'd_1'$.

\noindent {\bf Case 3:}  $d_1c_1'=0$, $d_1c_1'=0$, $d_1d_1'=0$.

\noindent {\bf Case 4:}  $d_1d_1=d_1'c_1$, $0= c_1d_1$, $d_1d_1=0$, $d_1c_1=0$.

\noindent {\bf Case 5:}   $0=c_1d_1$, $c_1d_1=0$, $0=d_1d_1$.

Those can be simplified to: $c_1d_1=c_1d_1'=c_1'd_1=d_1d_1=d_1d_1'=d_1'd_1'=0$.

$E_++E_-+c_1E+d_1VT=0$, $E_++E_-+c_1'E+d_1'S=0$. $R^{A_1}_3$ tells us that $(2+d_1)v_n+c_1v_{n+1}=0$.

\bigskip\noindent (4) For the type 2 invariant, if we ask $c_1=c_2=c_3=c_4=0$ and $\overline{c}_1=\overline{c}_2=\overline{c}_3=\overline{c}_4=0$, we get the followings.

\noindent {\bf Case 1:}  $d_1'd_2'=d_2'd_2'=0$.

\noindent {\bf Case 2:}  $d_1'd_2'+d_1'd_1'=0$, $d_2'd_2'+d_2'd_1'=0$.

\noindent {\bf Case 3:}  $d_1d_2'=0$, $d_2d_2'=d_1d_2'+d_1d_2'$, $d_1d_1'=0$, $0=d_2d_2'$.

\noindent {\bf Case 4:}  $d_2d_1+d_1d_1=0$,
$d_2d_2+d_1d_2=0$, $0=d_1d_2+d_1d_1$,

\noindent {\bf Case 5:}  $d_2d_2=d_1d_1$.

Those can be simplified.
 $d_2'd_2'=d_1'd_1'=d_1'd_2'=0$,
 $d_2d_2'=d_1d_1'=d_1d_2'=0$,
 $d_2d_2=d_1d_1=-d_1d_2$. They need to be complete by conjugation too. Check the following relation set.

\noindent $R^{A_2}_3$: $v_n=(c_1+c_2+c_3+c_4)v_{n+1}+(d_1+d_2)v_n$,

\noindent $v_n=(\overline{c}_1+\overline{c}_2+\overline{c}_3+\overline{c}_4)v_{n+1}+(\overline{d}_1+\overline{d}_2)v_n$,

\noindent $v_{n+1}=(\overline{c}_1+\overline{c}_2+\overline{d}_1+\overline{d}_2)v_{n}$.

Then we ask the following conditions hold. $v_{n+1}=v_n$ and $(d_1+d_2-1)v_n=0=(\overline{d}_1+\overline{d}_2-1)v_n$.

\bigskip There are other interesting invariants derived from type 2 invariants and their modified version. We don't discuss them here.

\begin{remark}
In general, we can regard the coefficient ring as commutative ring over rational numbers, then we can use tools like Gr$\ddot{o}$bner basis. Hence we can directly work with type one and type two invariants. However, the modified versions have infinitely many generators, so we need some technique to deal with them. We shall discuss them in the future.
\end{remark}

\section{An application}\label{sec:7}

Most knot invariants do not distinguish mutant knots. For example, the HOMFLY polynomial, the hyperbolic volume. It is not hard to show that neither type one nor type two invariant distinguishes mutant knots. However, in the following we shall show that the modified type one invariant can distinguish mutants. This invariant is very complicated, and is hard to compute by hand. We shall use a simplified version here. Hence the computation is much easier.

In the modified version of type one invariant, we let most variables to be zero, and only leave the followings. $E_++bE_-+d_2VT=0$, $E_++b'E_-+c_1'E=0$. Since we do not have other variables, we simply use $d$ to denote $d_2$, and $c$ to denote $c_1'$. The relations between the variables are: $b'c=bc,cd=0,dd=0,b^2d=d$. We shall modify it by writhe, so the Reidemeister invariance gives $\{h(w)+h(w-2)b+dh(w-1)\}v_n=0$.

 We simply ask $h(w)+h(w-2)b+dh(w-1)=0$ and $v_n$ commutes with other variables.

 If we let $x=h(1), y=h(2)$, the all other $h(w)$ are uniquely determined by the above equation. For example, $h(3)=-bx-dy$, $h(0)=-b^{-1}y-b^{-1}dx$.

 To calculate the invariant, recall that there are two functions $f(D)$ and $F(D)=h(w)f(D)$, where $F(D)$ is the knot invariant, and $f(D)$ satisfies the skein equations. For any monotone diagram $D$, the value for $f(D)$ is $h^{-1}(w)v_n$, and $F(D)=v_n$. Now for an arbitrary diagram $D$, if after some Reidemeister moves we get another diagram $D'$, then $f(D)=h^{-1}(w)F(D)=h^{-1}(w)F(D')=h^{-1}(w)h(w')f(D')$. We shall use this formula to simplify the calculation.

The famous first mutants pair in the knot table are the Conway knot $C$ ($11_{n34}$) of genus 3 and the Kinoshita¨CTerasaka ($11_{n42}$) knot KT of genus 2.

Calculation shows $F(KT)=h(-1)\{b^{-2}h^{-1}(3)+dh^{-1}(2)-b^{-1}d\}v_1$, and $F(C)=h(-1)\{bd+h^{-1}(-1)-dbh^{-1}(4)\}v_1$.

To compare those results, we need the Gr$\ddot{o}$bner basis. Let the commutative polynomial ring $R'=Q [ b,b',B,B',c,d ]$ and $R=Q  [ b, b',B,B',c,d,v_n, n=1,2,\cdots  ]= R'[v_n, n=1,2,\cdots  ]$. Here $Q$ denote the field of rational numbers. Let $I$ be the ideal in $R$ generated by $b'c-bc,cd,dd,b^2d-d,bB-1,b'B'-1$, and $I'$ be the ideal in $R'$ generated by $b'c-bc,cd,dd,b^2d-d,bB-1,b'B'-1$. Take the lex order $b'>b> B'>B>c>d>v_n>v_{n-1}$. Then a Gr$\ddot{o}$bner basis for $I'$ is $$G=\{ b'c-bc,cd,dd,b^2d-d,bB-1,b'B'-1, bd-dB,bcB'-c,B^2d-d, cB'-cB, b'cB-c\} .$$

Now let $A=R/I$, $A'=R'/I'$, $H=A[x,y], H'=A'[x,y]$. Let $S$ be the subset of non zero divisors in $A$, $T=S^{-1}A$. In $A$ and $A'$, we can regard $B$ as $b^{-1}$, $B'$ as $b'{}^{-1}$.

\begin{lem}
For any $f\in R$, $x+fd$ is not a zero divisor in $H$.
\end{lem}
\begin{proof}
Suppose that for some $e,g\in R$, $(x+fd)(e+gd)=0$. There are two cases. Case 1. $de=0$, then $x(e+gd)=0$. By definition of $H'$, $x$ is not a zero divisor. Hence $e+gd=0$.

Case 2. $de\neq 0$. Then $xe+d(fe+xg)=0$, multiply $d$ to both sides, we have $xde=0$, a contradiction.
\end{proof}

Using induction, we can prove that $h(2k+1)=(-b)^{k}x-k(-b)^{k-1}dy$, $h(2k+2)=(-b)^{k}y-k(-b)^{k-1}dx$, and $h(-2k)=(-b^{-1})^{k+1}y+(k+1)(-b^{-1})^{k+1}dx$,  $h(-2k-1)=(-b^{-1})^{k+1}x+(k+1)(-b^{-1})^{k+1}dy$ for $k>0$. For example, $h(3)=-bx-dy$, $h(4)=-by+bdx$, $h(0)=-b^{-1}y-b^{-1}dx$, $h(-1)=-b^{-1}x+b^{-2}dy=dy-b^{-1}x$.

It follows that all $h(k)$ are not zero divisors, and hence lie in $S$.

For any $u_1/s_1=u_2/s_2 \in T$, we have $s(s_1u_2-s_2u_1)=0$ for some $s \in S$. Since $S$ contains only non zero divisors, we have $s_1u_2-s_2u_1=0 \in H$.

Now we can compare $F(KT)$ and $F(C)$.

$F(KT)=h(-1)\{b^{-2}h^{-1}(3)+dh^{-1}(2)-b^{-1}d\}v_1 =h(-1)(b^{-2}\frac{1}{bx+dy}+d\frac{1}{y}-b^{-1}d)v_1$

$F(C)=h(-1)\{bd+h^{-1}(-1)-dbh^{-1}(4)\}v_1 = h(-1) (bd +\frac{1}{dy-b^{-1}x}-db\frac{1}{-by+bdx}) v_1$.

Suppose that $F(KT)=F(C)$, since $h(-1)$ is not a zero divisor, we shall have $b^{-2}\frac{1}{bx+dy}+d\frac{1}{y}-b^{-1}d =bd +\frac{1}{dy-b^{-1}x}-db\frac{1}{bdx-by}$.

Reduce the fractions to a common denominator and simplify, we shall have $2bdy^3-2dx^2y^2-(b^2-b^{-2})xy^2=0$ in $H$, which is not true. Hence $F(KT)\neq F(C)$. This invariant distinguishes mutants.

\bigskip There is a slightly weaker invariant can also distinguish mutants. This time we ask $h(n)=a^n$ for some invertible variable $a$.

We can now ask $b=b', B=B', c=d$, then we have $a^2+b+da=0$. Let $A>a>B>b>d$, then a Gr$\ddot{o}$bner basis is $b^2d-d,dd,bB-1,Bd-bd, aA-1, a^2+b+da, bA+a+d, A+aB+bd, dA+abd, bdA+ad$.

Now $a^3b^2F(KT)=a+da^2-bda^4=a-2bd$, $a^3b^2F(C)=bda^4-bd+a^3b^2=bd-ab^3$. Hence this invariant distinguishes mutants.

\section{Discussions}\label{sec:8}

The invariants are constructed in a similar pattern. Using the same pattern, we can construct more invariants.

$\\$
\noindent {\bf The Principle:} Whenever we get a choice for the skein relation, we calculate the equations $f_{pq}=f_{qp}$. Those equations defines an invariant at the diagram level. They almost give us a link invariant, and we only need to add a few equations to fit the Reidemeister invariance. There is usually a modified version and parameterized versions for the invariant.

$\\$
\noindent  With this in mind, we propose some more general link invariants.

One is constructed using the following skein relations. If the two strands are from same components, then

\noindent \resizebox{12.5cm}{!} {$E_{+}+bE_{-}+c_1E+\widehat{c}_1\overline{E}+c_2W+\widehat{c}_2\overline{W}+c_3HC+\widehat{c}_3\overline{HC}+c_4HT+\widehat{c}_4\overline{HT}+d_1VC+\widehat{d}_1\overline{VC}  +d_2VT+\widehat{d}_2\overline{VT}=0$}

\noindent If they are from different components, then
$$E_{+}+b'E_{-}+c_1'E+\widehat{c}_1'\overline{E}+c_2'W+\widehat{c}_2'\overline{W}+d'_1S+\widehat{d}'_1\overline{S}+d'_2N+\widehat{d}'_2\overline{N}=0$$

Here, $\widehat{E}$ is the same link diagram as $E$ with different orientations. For the link components containing the two strands, they have same orientation, but all other components are changed to opposite orientations. Similarly, one can define $\widehat{HC}$ and others. The proof is essentially the same, but the coefficient ring is far more complicated. We don't give it here.

$\\$
\noindent Another invariant has infinitely many variables. In stead of two or four skein relations, we use infinitely many skein relations. For example, if the two arrows/arcs in the local diagram are from the same link component, then

\noindent \resizebox{12.5cm}{!} {$E_{+}+b(cr, \mu)E_{-}+c_1(cr, \mu)E+c_2(cr, \mu)W+c_3(cr, \mu)HC+c_4(cr, \mu)HT+d_1(cr, \mu)VC+d_2(cr, \mu)VT=0.$}

\noindent If the two arrows/arcs are from different components, then
$$E_{+}+b'(cr, \mu)E_{-}+c_1'(cr, \mu)E+c_2'(cr, \mu)W+d'_1(cr, \mu)S+d'_2(cr, \mu)N=0.$$

Here $cr$ is the crossing number of diagram $E_+$, $\mu$ is the number of components of $E_+$. We call they the parameters. Here, for example, $b(1,1), b(1,2), b(2,1),b(3,2),\cdots $ are different independent variable. They can define a new knot invariant. Also, one can use other parameters, for example, the writhe $w$. The invariant then gives more direct information of crossing number and other knot invariants.

We call it the parametrization of type one invariant.

$\\$
\noindent We shall discuss the unoriented case in the next paper, and the generalized invariants in the future.

\begin{acknowledgements}
The author would like to thank Ruifeng Qiu, Jiajun Wang, Ying Zhang, Xuezhi Zhao, Hao Zheng, Teruhisa Kadokami for helpful discussions. The author also thanks all the organizers of 2009 summer school of knot theory in ICTP, Italy. He really had a great time there and got many inspirations. The author also thanks NSF of China for providing traveling fee for the trip, and all the efforts made by ICTP, Italy.
\end{acknowledgements}

\nocite{*}

\bibliographystyle{spmpsci}      
\bibliography{ref}

\end{document}